\documentclass[11pt]{amsart}
\usepackage[english]{babel}
\usepackage{graphicx}
\newtheorem{theorem}{Theorem}[section]

\theoremstyle{definition}

\theoremstyle{remark}
\newtheorem{remark}{Remark}
\newtheorem{example}{Example}
\numberwithin{equation}{section}
\newcommand{\real}{\mathbb R}
\newcommand{\R}{\mathbb R}

\def\natu{\mathbb N}
\def\N{\mathbb N}
\def\dis{\displaystyle}
\def\sobol{H^{1}(0,T)}
\def\intl{\int_{0}^{T}}
\def\landan{\Lambda_n}
\def\landa0{\Lambda_0}
\def\landaene{\lambda_{2n-1}}
\def\periodica{L_{T}(\real,\real)}
\def\periodicapi{L_{\pi}(\real,\real)}
\def\landaenea{\lambda_n (a)}
\def\betaant{\beta_p ^{ant}}
\def\comega{\overline{\Omega}}
\def\intomega{\int_{\Omega}}

\def\landan{\Lambda_n}

\newcommand{\norm}[1]{\left\Vert#1\right\Vert}

\def\afirm#1#2#3{\skipaline\noindent{\bf #1} \hfill
\begin{tabular}{p{#2}}{\sl #3} \end{tabular}\hfill $ $\skipaline}
\def\skipaline{\vskip12pt plus 1pt}

\begin{document}

\title[A mathematical excursion through Lyapunov inequalities] {An applied mathematical excursion through Lyapunov
inequalities, classical analysis and differential equations}
\author{Antonio Ca\~{n}ada}
\author{Salvador Villegas}
\thanks{The authors have been supported by the Ministry of Education and
Science of Spain (MTM2008.00988) and by Junta de Andalucia
(FQM116)}
\address{Departamento de An\'{a}lisis
Matem\'{a}tico, Universidad de Granada, 18071 Granada, Spain.}
\email{acanada@ugr.es, svillega@ugr.es}


%\date{}%
%\commby{}%
% ----------------------------------------------------------------
\begin{abstract}
Several different problems make the study of the so called
Lyapunov type inequalities of great interest, both in pure and
applied mathematics. Although the original historical motivation
was the study of the stability properties of the Hill equation
(which applies to many problems in physics and engineering), other
questions that arise in systems at resonance, crystallography,
isoperimetric problems, Rayleigh type quotients, etc. lead to the
study of $L_p$ Lyapunov inequalities ($1\leq p\leq \infty$) for
differential equations. In this work we review some recent results
on these kinds of questions which can be formulated as optimal
control problems. In the case of Ordinary Differential Equations,
we consider periodic and antiperiodic boundary conditions at
higher eigenvalues and by using a more accurate version of the
Sturm separation theory, an explicit optimal result is obtained.
Then, we establish Lyapunov inequalities for systems of equations.
To this respect, a key point is the characterization of the best
$L^p$ Lyapunov constant for the scalar given problem, as a minimum
of some especial (constrained or unconstrained) variational
problems defined in appropriate subsets of the usual Sobolev
spaces. For Partial Differential Equations on a domain $\Omega
\subset \real^N$, it is proved that the relation between the
quantities $p$ and $N/2$ plays a crucial role in order to obtain
nontrivial $L_p$ Lyapunov type inequalities (which are called
Sobolev inequalities by many authors). This fact shows a deep
difference with respect to the ordinary case. Combining the linear
results with Schauder fixed point theorem, we can obtain some new
results about the existence and uniqueness of solutions for
resonant nonlinear problems for ODE or PDE, both in the scalar
case and in the case of systems of equations.
\end{abstract}

\maketitle

\section{Introduction}
The Hill equation
\begin{equation}\label{neq1}
u''(t) + a(t)u(t) = 0, \ t \in \R
\end{equation}
where the function $a$ satisfies
\begin{equation}\label{nhipotesisa}
a: \R \rightarrow \R \ \ \mbox{is}\ T-\mbox{periodic and } \ a \in
L^1 (0,T),
\end{equation}
models many phenomena in applied sciences (\cite{cole},
\cite{hale}, \cite{magnus}). In a broader sense, Hill's equation
connotes the class of homogeneous, linear, second order
differential equations with real and periodic coefficients. In the
preface to the book on Hill's equation, by Magnus and Winkler
\cite{magnus}, we can read:

$ $

{\it ``There exist hundreds of applications of Hill's equation to
problems in engineering and physics, including problems in
mechanics, astronomy, the theory of electric circuits, of the
electric conductivity of metals, of the cyclotron and new
applications are continually discovered''.}

$ $

In what follows, we denote by $L_T (\R,\R)$ the set of functions
$a(\cdot)$ satisfying (\ref{nhipotesisa}).

$ $

The study of stability properties of (\ref{neq1}) is of special
interest. Whenever all solutions of (\ref{neq1}) are bounded, we
say that (\ref{neq1}) is stable; otherwise we say that it is
unstable. Floquet theory assures that such stability properties
depend on characteristic multipliers, which are defined using any
fundamental matrix of the given equation. This important
theoretical result requires the knowledge of all solutions of
(\ref{neq1}).

In early twentieth century Lyapunov proved (\cite{lyap},
\cite{magnus}) that if
\begin{equation}\label{nhill2}
0 \prec a, \ \intl a (t) \ dt \leq \frac{4}{T},
\end{equation}
then (\ref{neq1}) is stable.  Here, for $c,d \in L^1 (0,T),$ we
write $c \prec d$ if $c(t) \leq d(t)$ for a.e. $t \in [0,T]$ and
$c(t) < d(t)$ on a set of positive measure.

Lyapunov's result is remarkable, among others, for two main
reasons: first, one can check (\ref{nhill2}) directly from the
equation (\ref{neq1}) and second, it is optimal in the following
sense (\cite{magnus}): for any $\varepsilon \in \R^+$, there are
unstable differential equations (\ref{neq1}) with $a$ satisfying
(\ref{nhipotesisa}), for which
$$
0 \prec a,  \ \intl a(t) \ dt \leq \frac{4}{T} + \varepsilon
$$
Condition (\ref{nhill2}) has been generalized in several ways
(\cite{borg}, \cite{krein}). More recently the authors provide in
\cite{zhangli} optimal stability criteria by using $L^{p}$ norms
of $a^+,$ $1 \leq p \leq \infty$.

$ $

The parametric equation
\begin{equation}\label{n2hilleq}
u''(t) + (\mu+a(t)) u(t) = 0, \ \mu \in \real,
\end{equation}
plays a fundamental role in the study of the stability properties
of the Hill equation (\ref{neq1}). In fact, if $\lambda_i (a), \ i
\in \N \cup \{ 0 \}$ and $\tilde{\lambda}_i (a), \ i \in \N,$
denote, respectively, the eigenvalues of (\ref{n2hilleq}) for the
periodic ($u(0)-u(T) = u'(0)-u'(T) = 0$) and antiperiodic
($u(0)+u(T) = u'(0)+u'(T) = 0$) problem, then Lyapunov and Haupt
proved (see (\cite{cole}, \cite{hale} and \cite{magnus} for
historical and mathematical details) that
\begin{equation}\label{n2004097}
\lambda_0 (a) < \tilde{\lambda}_1 (a) \leq \tilde{\lambda}_2 (a) <
\lambda_{1}(a) \leq \lambda_{2}(a) < \tilde{\lambda}_{3}(a) \leq
\tilde{\lambda}_{4}(a) < \lambda_3 (a) \leq \ldots
\end{equation}
and that equation (\ref{n2hilleq}) is stable if
\begin{equation}\label{n2004095}
\mu \in (\lambda_{2n}(a),\tilde{\lambda}_{2n+1}(a)) \cup
(\tilde{\lambda}_{2n+2}(a),\lambda_{2n+1}(a))
\end{equation}
for some $n \in \N \cup \{ 0 \}.$ As a consequence, if either
\begin{equation}\label{1403111}
\lambda_{2n}(a) < 0 < \tilde{\lambda}_{2n+1}(a) \end{equation} or
\begin{equation}\label{1403112}
\tilde{\lambda}_{2n+2}(a) < 0 < \lambda_{2n+1}(a) \end{equation}
then (\ref{neq1}) is stable. In this case, we say that $\mu = 0$
belongs to the $n^{th}$ stability zone of (\ref{n2hilleq}). In
particular, the Lyapunov's conditions (\ref{nhill2}) imply
$\lambda_0 (a) < 0 < \tilde{\lambda}_1 (a).$

$ $

Except in very special cases, it is not an easy task to obtain
some information on the sign of the previous eigenvalues. One
possibility is to use the variational characterization, but in
general it is a difficult question, especially at higher
eigenvalues.

$ $

It is at this point where the so called ``Lyapunov inequalities''
can play an important role. In fact, the problem of the
determination of sufficient conditions to have the desired
properties (\ref{1403111}) or (\ref{1403112}), is intimately
connected with the problem of finding  the so called best Lyapunov
constant for certain linear and homogeneous boundary value
problems (Theorem \ref{pt1} and Theorem \ref{antipt1} in the
second section). Illustrate this with an example for periodic
boundary conditions.

Let us consider the problem
\begin{equation}\label{np1}
u''(t) + a(t)u(t) = 0, \ t \in (0,T), \ u(0)-u(T) = u'(0)- u'(T) =
0
\end{equation}
where $a\in \periodica$. We denote
$$ \lambda_0 (0) = \lambda_0 = 0, \ \lambda_{2n-1}(0) =
\lambda_{2n}(0)= \lambda_{2n-1} = \lambda_{2n}=  (2n)^2 \pi^2/T^2,
\ n \in \N.$$ Then, if $n \in \natu$ is fixed, we introduce the
set $\landan$ as
\begin{equation}\label{1503111}
\landan = \{ a \in \periodica: \lambda_{2n-1} \prec a \ \mbox{and}
\ (\ref{np1}) \ \mbox{has nontrivial solutions} \ \}
\end{equation}
and define the Lyapunov constant (at higher eigenvalues) as
\begin{equation}
\gamma_{1,n} = \inf_{a \in \landan } \ \Vert a  \Vert_{L^1 (0,T)}
\end{equation}
Then if a given function $a \in \periodica$ satisfies
$$
\lambda_{2n-1} \prec a, \ \Vert a \Vert_{L^1(0,T)} < \gamma_{1,n}
$$
we have
$$
\lambda_{2n}(a) < 0 < \lambda_{2n+1}(a).
$$
In fact, since $\lambda_{2n-1} \prec a,$ by using the variational
characterization of the periodic eigenvalues, we trivially have
$\lambda_{2n}(a) < \lambda_{2n}(\lambda_{2n-1}) = 0.$ To prove
that $\lambda_{2n+1} (a) >0$ we can use the following continuation
method (which, to the best of our knowledge, has been introduced
by the authors in \cite{cavijmaa}): let us define the continuous
function $g: [0,1] \rightarrow \real$ by
$$
g(\varepsilon) = \lambda_{2n+1}(a_{\varepsilon}(\cdot))
$$
where $a_{\varepsilon}(x) = \lambda_{2n-1} + \varepsilon (a(\cdot)
- \lambda_{2n-1}).$ Then $g(0) =\lambda_{2n+1}\left (
\lambda_{2n-1} \right )= \lambda_{2n+1} - \lambda_{2n-1} > 0.$
Moreover, for each $\varepsilon \in (0,1),$ the function
$a_{\varepsilon}$ satisfies: $a_{\varepsilon} \in \periodica, \
\lambda_{2n-1} \prec a_{\varepsilon}$ and $\Vert a_{\varepsilon}
\Vert_{L^1 (0,T)} < \gamma_{1,n}.$ Therefore, taking into account
the definition of $\gamma_{1,n},$ the real number zero is not an
eigenvalue of the eigenvalue problem
$$
u''(t) + (\mu + a_{\varepsilon} (t))u(t) = 0, \ t \in (0,T), \ \
u(0)-u(T) = u'(0) - u'(T) = 0.
$$
In particular, $\lambda_{2n+1} (a_{\varepsilon}(\cdot)) \neq 0$,
i.e. $g(\varepsilon) \neq 0, \ \forall \ \varepsilon \in (0,1]$
and consequently, we obtain $g(1) = \lambda_{2n+1}(a)> 0.$

These arguments can be used not only for $L_1-$norms but also for
$L_p$-norms with $1 \leq p \leq \infty$ and for different boundary
conditions such as antiperiodic, Dirichlet, Neumann, etc. As we
show in Theorem \ref{pt1} and Theorem \ref{antipt1}, we use a
common procedure to study this kind of problems and also (see the
previous commentaries to Theorem \ref{pt1}) for obtaining possible
applications, which include not only stability theory, but also
resonant (linear or nonlinear) boundary value problems (see
Theorem \ref{nolineal} and Theorem \ref{nolinealp}). Of course,
each boundary value problem has its peculiarities and, moreover,
the cases $p=1,$ $p= \infty$ and $1 < p < \infty$ are generally
different from a qualitative point of view, specially in the
problems of the calculus of variations associated to them (see
Theorem \ref{l2}). These problems are, among other things, of
interest in crystallography (\cite{dacoro2}, \cite{dacoro1}).

$ $

One last general comment before describing the contents of the
sections of the paper: this type of question can be seen as an
optimal control problem and some authors have used the Pontryagin
maximum principle in its study (\cite{giau},
\cite{huaizhongyongjde94}, \cite{liwang}, \cite{wanglina95}). For
example, for the periodic problem, the admissible control set is
$\landan$ defined in (\ref{1503111}) and the functional that we
want to minimize is $\Phi: \landan \rightarrow \R,$ defined as
$\Phi (a) = \intl a(t) \ dt, \ \forall \ a \in \landan.$ However,
we caution that the condition
\begin{equation}
(\ref{np1}) \ \mbox{has nontrivial solutions}
\end{equation}
is difficult to handle from a mathematical point of view and this
is the main difficulty of the problem. We will overcome these
difficulties by using different methods of classical analysis such
as a modified version of Sturm separation theorem, constrained or
unconstrained minimization problems of the calculus of variations
(depending on the boundary value problem that we are considering),
Euler equation, Lagrange multiplier theorem, continuation methods,
Schauder fixed point theorem, etc.

$ $

In the second section of this paper we consider scalar ordinary
differential equations. First we study the best $L_1$ Lyapunov
constant at higher eigenvalues for the periodic and the
antiperiodic boundary value problem. Then, we apply these results
to obtain some new sufficient conditions on the $n^{th}$ stability
zone of the Hill equation (\ref{neq1}).

$ $

The third section is dedicated to systems of equations like
\begin{equation}\label{n2ieq35}
u''(t) + P(t)u(t)= 0,\ t \in \real,
\end{equation}
where the matrix function $P(\cdot)$ is $T-$periodic. Despite its
undoubted interest (see references \cite{cole}, \cite{hale},
\cite{krein}, \cite{magnus} to get an idea of the many
applications of this kind of systems), there are not many studies
on this subject. We must remark that the results we show in the
scalar case are mainly based on an exact knowledge about the
number and distribution of the zeros of the corresponding
solutions, but to the best of our knowledge, for systems of
equations like (\ref{n2ieq35}), we do not know  similar results to
those given in Theorem \ref{pt1} and Theorem \ref{antipt1}.

To extend the obtained results in the scalar case to systems of
equations, a key point is the characterization of the best $L_p$
Lyapunov constant for the scalar given problem, as a minimum of
some especial minimization problems defined in appropriate subsets
$X_p, \ 1 \leq p \leq \infty,$ of the Sobolev space $H^1 (0,T). $
For the Dirichlet problem this was done by Talenti
(\cite{talenti}) and for the Neumann and periodic problem this was
done by the authors in \cite{camovimia} and \cite{cavina}(see also
\cite{dacoro2} for a different treatment and motivation of the
problem). Since $0$ is the first eigenvalue for Neumann and
periodic boundary problems, it is necessary to impose an
additional restriction to the definition of the spaces $X_p, \ 1
\leq p \leq \infty$ in the case of Neumann and periodic conditions
(see Remark \ref{080311} below), and therefore we will have
constrained minimization problems. This is not necessary in the
case of Dirichlet or antiperiodic boundary ones where we will find
unconstrained minimization problems.

In regard to stability properties of (\ref{n2ieq35}), which will
be considered also in Section 3, the results proved by Krein in
\cite{krein} show that the problem is closely related to Lyapunov
inequalities. In fact, under appropriate hypotheses on the
function matrix $P(\cdot)$ (see Subsection 3.2), the stability
properties of (\ref{n2ieq35}) strongly depend on the fact that the
smallest positive eigenvalue of the antiperiodic eigenvalue
problem
\begin{equation}\label{n2eq37}
u''(t) + \lambda P(t)u(t)= 0, \ t \in \real, \ u(0)+u(T) =
u'(0)+u'(T) = 0
\end{equation}
be greater than one. We establish some new conditions to get this
last property.

$ $

Section 4 will be devoted to the study of $L_p$ Lyapunov
inequalities and some applications to nonlinear resonant problems
for scalar partial differential equations. More precisely, we
consider the linear problem with Neumann boundary conditions
\begin{equation}\label{nceq6}
\left.
\begin{array}{cl}
\Delta u(x)+ a(x)u(x)= 0 & x\in \Omega \, \\
 \frac{\partial u}{\partial n}(x)=0\, \ \ \ \ \ \ &  x\in \partial \Omega
\end{array}
\right\}
\end{equation}
\noindent where $\Omega \subset \real^N$ ($N\geq 2$) is a bounded
and regular domain and  $\dis \frac{\partial}{\partial n}$ is the
outer normal derivative on $\partial \Omega.$

According to our opinion, the work by the authors
\cite{camovibrezis} and \cite{camovijfa} showed for the first time
significant differences from the ordinary case. The most important
is, without doubt, the following fact: if $N=2,$ there is $L_p$
Lyapunov inequality if and only if $1 < p \leq \infty$ and if $N
\geq 3,$ there is $L_p$ Lyapunov inequality if and only if $N/2
\leq p \leq \infty.$ In particular, we must remark that there is
no $L_1$ Lyapunov inequality in the PDE case and that the best
constant is attained if $N/2 < p \leq \infty.$

A pioneering work for the case of Dirichlet boundary conditions
$u(x) = 0, \ x \in \partial \Omega,$ was undoubtedly the Talenti
article \cite{talenti}, where the author relates this kind of
inequalities (called Sobolev Inequalities in this and in other
many papers in PDE), with the classic isoperimetric problem for
subsets of Euclidean space. Talenti made a detailed study of the
case where $\Omega$ is substituted by the whole space $\R^N$, $N
\geq 3$ and $p = N/2$ (the critical value), including the study of
some additional symmetric properties for those functions where the
considered infimum is attained.

By using our methods, in the case of Dirichlet conditions it is
possible to obtain analogous results and if $N/2<p<\infty$, the
Lyapunov constant for the Dirichlet problem, $\beta^D_{p},$ can be
characterized variationally as
$$\beta^D_p=\inf_{u\in H_0^1(\Omega)\setminus\{ 0\} }\frac{\dis \int_\Omega \vert \nabla
u\vert^2}{\left ( \dis \int_{\Omega} \vert u \vert
^{\frac{2p}{p-1}} \right )^{\frac{p-1}{p}}}$$

$ $

\noindent If $\Omega$ is, moreover, a radial domain, previous
minimization problem is related to a more general one which
involves Rayleigh quotient
$$
\frac{\dis \int_\Omega \vert \nabla u\vert^2}{\left ( \dis
\int_{\Omega} \rho (x) \vert u \vert ^{\frac{2p}{p-1}} \right
)^{\frac{p-1}{p}}}
$$
where $\rho \in L^{q}(\Omega), \ q = N(p-1)/(2p-N),$ is a positive
function. This has been used in the study of the existence of
nonsymmetric ground states of symmetric problems for nonlinear
PDE's (see \cite{br}, \cite{brni} and \cite{smwisu}).

In Section 4, first we treat the scalar linear problem
(\ref{nceq6}). Then we apply the linear study to nonlinear
resonant problems by using the Schauder fixed point theorem (see
Theorem \ref{nolineal}). Obviously, similar results may be proved
for nonlinear and resonant ordinary differential equations for
each value $1 \leq p \leq \infty.$

$ $

In Section 5,  we deal with elliptic systems of the form

\begin{equation}\label{ns1p}
\Delta u(x)+ A(x)u(x) = 0,  \ x\in \Omega,  \ \ \frac{\partial
u(x)}{\partial n}=0, \  x\in \partial \Omega,
\end{equation}

\noindent where $A(\cdot)$ is a real $n\times n$ symmetric and
continuous matrix satisfying some additional restrictions. To our
knowledge, there are no previous work on $L_p$ Lyapunov
inequalities for elliptic systems like (\ref{ns1p}) if $p \neq
\infty$ (see \cite{ba} and \cite{kawa}, Section 5, for the case $p
= \infty$).

As in the ordinary case, to extend to systems of equations the
obtained results in the scalar case, a key point is the
characterization of the best $L_p$ Lyapunov constant for the
scalar given problem, as a minimum of some especial minimization
problems defined in appropriate subsets $X_p, \ N/2 <  p \leq
\infty,$ of the Sobolev space $H^1 (\Omega). $

By combining the linear study with Schauder fixed point theorem,
we show some applications to nonlinear resonant problems. In
particular, and for Neumann boundary conditions, we obtain a
generalization for systems of equations of the main result given
in \cite{mawawi} where the authors treated the scalar case and
where they use in the proof the duality method of Clarke and
Ekeland.

$ $

Finally, we present some remarks on other interesting aspects not
treated here and also on some open problems.

\section{Scalar ordinary differential equations}
This section deals with $L_1$ Lyapunov inequality for the periodic
and the antiperiodic boundary value problem at higher eigenvalues.
As an application, we obtain some new stability properties on the
Hill equation (\ref{neq1}).

We explain in more detail what we mean by $L_1$ Lyapunov
inequality at higher eigenvalues. Remember that
$L_{T}(\real,\real)$ denotes the set of $T-$periodic functions $a:
\R \rightarrow \R,$ such that $a|_{[0,T]} \in L^1 (0,T).$ Then, we
will consider the periodic boundary value problem
\begin{equation}\label{p1}
u''(t) + a(t)u(t) = 0, \ t \in (0,T), \ u(0)-u(T) = u'(0)- u'(T) =
0
\end{equation}
It is very well known (\cite{cole}, \cite{hale}, \cite{hart}) that
the eigenvalues for
\begin{equation}\label{p2bisbis}
 u''(t) + (\lambda+a(t)) u(t) = 0, \ t \in (0,T), \ u(0)-u(T) = u'(0)- u'(T) =
0 \end{equation} form a sequence $\landaenea, \ n \in \natu \cup
\{ 0 \},$ such that
\begin{equation}\label{1604091}
\lambda_0 (a) < \lambda_1 (a) \leq \lambda_2 (a)< \ldots <
\lambda_{2n-1}(a) \leq \lambda_{2n}(a) < \ldots
\end{equation}
with $\lambda_0(a)$ simple and such that if $\phi_n$ is the
corresponding eigenfunction to $\lambda_n(a),$ then $\phi_0$ has
no zeros in $[0,T]$ and $\phi_{2n-1}$ and $\phi_{2n}$ have exactly
$2n$ zeros in $[0,T)$. In particular, $\lambda_0 = \lambda_0 (0) =
0, \ \lambda_{2n-1} = \lambda_{2n}= \lambda_{2n-1}(0) =
\lambda_{2n}(0)= (2n)^2 \pi^2/T^2, \ n \in \natu.$

Similarly, we will consider the antiperiodic boundary value
problem
\begin{equation}\label{ap1}
u''(t) + a(t)u(t) = 0, \ t \in (0,T), \ u(0)+u(T) = u'(0)+ u'(T) =
0
\end{equation}
where $a \in \periodica.$ It is also very well known (\cite{cole},
\cite{hale}, \cite{hart}) that the eigenvalues for
\begin{equation}\label{antip2bisbis}
 u''(t) + (\tilde{\lambda}+a(t)) u(t) = 0, \ t \in (0,T), \ u(0)+u(T) = u'(0)+ u'(T) =
0 \end{equation} form a sequence $\tilde{\lambda}_{n} (a), \ n \in
\natu,$ such that
\begin{equation}\label{anti1604091}
\tilde{\lambda}_1 (a) \leq \tilde{\lambda}_2 (a)< \ldots <
\tilde{\lambda}_{2n-1}(a) \leq \tilde{\lambda}_{2n}(a) < \ldots
\end{equation}
and if $\tilde{\phi}_n$ is the corresponding eigenfunction to
$\tilde{\lambda}_n(a),$ then $\tilde{\phi}_{2n-1}$ and
$\tilde{\phi}_{2n}$ have exactly $2n-1$ zeros in $[0,T)$. In
particular, $\tilde{\lambda}_{2n-1}(0) = \tilde{\lambda}_{2n}(0)=
(2n-1)^2 \pi^2/T^2, \ n \in \natu.$ We will denote
$\tilde{\lambda}_i = \tilde{\lambda}_i (0), \ \forall \ i \in
\natu.$

$ $

It is clear that if the function $a(\cdot)$ is not a constant, it
is not possible, in general, to calculate explicitly the mentioned
eigenvalues. Specifically, it is not easy to establish that the
real number $0$ is not an eigenvalue of (\ref{p2bisbis}) (or
(\ref{antip2bisbis})). Perhaps, for the periodic problem
(\ref{p1}), the following result is intuitive and understandable:
if there exists some $n \in \N$ such that $\lambda_{2n-1} \prec a
\prec \lambda_{2n+1}$ (or $\lambda_{0} \prec a \prec
\lambda_{1}$), then the unique solution of (\ref{p1}) is the
trivial one, i.e., $0$ is not an eigenvalue of (\ref{p2bisbis})
(see \cite{laze}). A similar result is true for the antiperiodic
problem.

$ $

The previous criterion may be seen as a $L_\infty-L_\infty$
criterion in the sense that the function $a$ is bounded, both
above and below, by consecutive different eigenvalues
$\lambda_{2n-1}, \lambda_{2n+1}.$ It allows a weak interaction
between the function $a$ and the corresponding spectrum. Our next
purpose is to show a $L_\infty-L_1$ criterion. It will allow the
function $a$ to cross an arbitrary number of eigenvalues as long
as certain $L_1$-norms are controlled.

$ $

More specifically, if $n \in \natu $ is fixed, we introduce for
the periodic problem, the set $\landan$ as
\begin{equation}\label{p2806073}
\landan = \{ a \in \periodica : \lambda_{2n-1} \prec a \
\mbox{and} \ (\ref{p1}) \ \mbox{has nontrivial solutions} \ \}
\end{equation}
Our main interest is in having a broader knowledge as possible
about the constant
\begin{equation}\label{pconstante}
\gamma_{1,n} = \inf_{a \in \landan } \ \Vert a  \Vert_{L^1 (0,T)}
\end{equation}

Analogously, if $n \in \natu$ is fixed, we can introduce for the
antiperiodic problem the set $\tilde{\Lambda}_{n}$ as
\begin{equation}
\tilde{\Lambda}_{n} = \{ a \in \periodica: \tilde{\lambda}_{2n-1}
\prec a \ \mbox{and} \ (\ref{ap1}) \ \mbox{has nontrivial
solutions} \ \}
\end{equation}
and, as for the periodic problem, our main interest is in having a
broader knowledge as possible about the constant
\begin{equation}\label{anticonstante}
\tilde{\gamma}_{1,n} \equiv \inf_{a \in \tilde{\Lambda}_{n} } \
\Vert a  \Vert_{L^1 (0,T)}
\end{equation}

$ $

The above ideas clarify the meaning of the expression ``$L_1$
Lyapunov inequalities at higher eigenvalues''.

Additionally, we prove some qualitative properties on the sign of
the eigenvalues $\lambda_{2n}(a), \lambda_{2n+1}(a),
\tilde{\lambda}_{2n}(a), \tilde{\lambda}_{2n+1}(a),$ which will be
very important in the applications to resonant problems and in the
study of stability properties of linear periodic equations. To
prove these qualitative properties, we use a curious continuation
method applied to an especial continuous function which does not
change its sign in its domain (an interval).

\subsection{Lyapunov constants for the periodic and the antiperiodic boundary value problem}

First, we deal with the constant $\gamma_{1,n}$ defined in
(\ref{pconstante}). Note that if $a \in \landan,$ and $u$ is any
nontrivial solution of (\ref{p1}), then $u$ is not a constant
function and $u$ must have a zero in the interval $[0,T].$  If $r
\in [0,T]$ is such that $u(r) =0,$ the periodic and nontrivial
function $v(t) = u(r+t)$ satisfies $v(0) = v(T) = 0,$ $v''(t) +
a(r +t)v(t) = 0, \ t \in (0,T)$ and $\Vert a(r + \cdot) \Vert_{L^1
(0,T)}= \Vert a(\cdot) \Vert_{L^1 (0,T)}.$ This observation does
possible to consider in our study only those situations where $u$
is such that $u(0) = u(T) = 0.$ Since $a \in \landan,\ n \in
\natu,$ it is clear that between two consecutive zeros of the
function $u$ there must exists a zero of the function $u'$ and
between two consecutive zeros of the function $u'$ there must
exists a zero of the function $u.$

$ $

In the next theorem, we provide a comprehensive review of the
control problem which we have described in (\ref{p2806073}) and
(\ref{pconstante}) and we show some qualitative properties of the
eigenvalues $\lambda_{2n}(a), \lambda_{2n+1}(a)$ associated to
those functions $a \in \periodica$ such that $\Phi (a) = \intl
a(t) \ dt \leq \gamma_{1,n}.$ We consider these qualitative
properties (and similar ones which will be obtained for the
antiperiodic boundary value problem) interesting in themselves.

$ $

In the proof, the following fundamental ideas are used:
\begin{enumerate}
\item A careful analysis on the number and distribution of zeros
of nontrivial solutions $u$ of the equation (\ref{p1}) (or the
equation (\ref{ap1})) and their first derivatives $u',$ including
an optimal estimation about the corresponding distances between
the zeros of $u$ and $u'.$ To this respect, we compare our problem
with other one with mixed boundary conditions. This allows us to
obtain a more precise information than that obtained when the
classical Sturm separation theorem is used, where the problem is
compared with the case of Dirichlet boundary conditions
(\cite{hart}) and therefore, the obtained information is only on
the zeros of $u$.

\item The previous step will motivate the study of the following
minimization problem ($H^1 (a,b)$ denotes the usual Sobolev
space):

Let $H = \{ u \in H^1 (a,b) : u(a) = 0, u(b) \neq 0 \}$ and $J: H
\rightarrow \R$ be defined by
\begin{equation}\label{02073}
J(u) = \frac{\int_{a}^{b} u'^2 - M \int_{a}^{b} u^2 }{u^2(b)}
\end{equation}
where $a < b$ and $0 < M \leq \frac{\pi^2}{4(b-a)^2},$ the first
(or principal) eigenvalue of the eigenvalue problem with mixed
boundary conditions
\begin{equation}\label{bis2806074} u''(t) + \mu u(t) = 0, \ t \in
(a,b), \ u(a) = 0, \ u'(b) = 0. \end{equation} Then, with the help
of classical and standard methods of the calculus of variations
(existence of minimizing sequences and the study of the
Euler-Lagrange equation), we obtain that $c \equiv \inf_{u \in H}
\ J(u),$ is attained and
\begin{equation}\label{050420111}
c = M^{1/2}\cot (M^{1/2}(b-a))
\end{equation}
Moreover, if $u \in H,$ then $J(u) = c \Longleftrightarrow u(t) =
\displaystyle  k \sin (M^{1/2}(t-a))$ for some non zero constant
$k.$

\item The final step is the study of a simple minimization problem
which is given by a finite sum of functionals of the previous
type. More precisely, given any $r\in\N $ and $S\in\R^+$
satisfying $r \pi > 2S$, let
$$
Z = \{ z= (z_0,z_1,\ldots,z_{r-1}) \in (0,\pi/2]^r: \
\sum_{i=0}^{r-1} z_i = S \}
$$
If $F: Z \rightarrow \real$ is defined by
\begin{equation}\label{1603111}
F(z) = \sum_{i=0}^{r-1} \ \cot \ z_i,
\end{equation}
we will show that $ \displaystyle \inf_{z \in Z} \ F(z)$ is
attained and its value is $ r \cot \frac{S}{r}.$ Moreover, $z \in
Z$ is a minimizer if and only if $z_i=\frac{S}{r}, \ \forall \ 0
\leq i \leq r-1.$ To obtain these conclusions, we basically use
Lagrange multiplier theorem.

\end{enumerate}

\begin{theorem}\label{pt1}
Let $n \in \natu$ and $a \in \landan$ be given and $u$ any
nontrivial solution of (\ref{p1}) such that $u(0) = u(T) = 0.$ If
the zeros of $u$ in $[0,T]$ are denoted by $0 = t_0 < t_2 < \ldots
< t_{2m} = T$ and the zeros of $u'$ in $(0,T)$ are denoted by $
t_1 < t_3 < \ldots < t_{2m-1}, $ then:
\begin{enumerate}
\item $t_{i+1} - t_{i} \leq \frac{T}{4n}, \ \forall \ i: \ 0 \leq
i \leq 2m-1.$ Moreover, at least one of these inequalities is
strict. \item $m$ is an even number and $m \geq 2(n+1).$  Any even
value $m\geq 2(n+1)$ is possible. %
%\item
%\begin{equation}\label{02071}
%\Vert a - \landaene \Vert_{L^1 (x_i,x_{i+1})} \geq \frac{2n\pi}{T}
%\cot (\frac{2n\pi}{T}(x_{i+1}-x_i)),\ 0 \leq i \leq 2m-1.
%\end{equation}
\item $\gamma_{1,n}$ is not attained and
\begin{equation}\label{p1009071}
\gamma_{1,n}  = T \lambda_{2n-1} + \frac{8\pi n (n+1)}{T} \cot
\frac{n\pi}{2(n+1)}
\end{equation}
\item If $a \in \periodica$ satisfies
\begin{equation}\label{positivo}
\lambda_{2n-1} \prec a, \ \Vert a \Vert_{L^1(0,T)} \leq
\gamma_{1,n} \end{equation} then
\begin{equation}\label{1604096}
\lambda_{2n}(a) < 0 < \lambda_{2n+1}(a)
\end{equation}
\end{enumerate}
\end{theorem}

\begin{proof} We only summary the main ideas (see \cite{cavijems}, Theorem 2.1 and \cite{cavijmaa}, Theorem 2.1
for the details).

$ $

First, let $i, \ 0 \leq i \leq 2m-1,$ be given. Then, function $u$
satisfies either the problem \begin{equation}\label{2806076}
u''(t) + a(t)u(t) = 0, \ t \in (t_i,t_{i+1}), \ \ u(t_i) = 0, \
u'(t_{i+1}) = 0, \end{equation} or the problem
\begin{equation}\label{2806077} u''(t) + a(t)u(t) = 0, \ t \in (t_i,t_{i+1}),
\ \ u'(t_i) = 0, \ u(t_{i+1}) = 0. \end{equation} Let us assume
the first case and let us denote by $\mu_1 ^i =
\frac{\pi^2}{4(t_{i+1}-t_i)^2}$ and $\varphi_1^i (t) = \sin
\frac{\pi (t-t_i)}{2(t_{i+1}-t_i)},$ respectively, the principal
eigenvalue and eigenfunction of the eigenvalue problem
\begin{equation}\label{2806074} v''(t) + \mu v(t) = 0, \ t \in
(t_i,t_{i+1}), \ v(t_i) = 0, \ v'(t_{i+1}) = 0. \end{equation}

Choosing $\varphi_1 ^i$ as test function in the weak formulation
of (\ref{2806076}) and $u$ as test function in the weak
formulation of (\ref{2806074}) for $\mu = \mu_1 ^i$ and $v
=\varphi_1^i,$ we obtain
\begin{equation}\label{2806078}
\int_{t_i}^{t_{i+1}} (a(t) - \mu_1^i)u(t) \varphi_1^i (t) \ dt =
0.
\end{equation}
This last relation and the hypothesis $\lambda_{2n-1} \prec a$,
prove the first part of the theorem and also the relation $m \geq
2n+1$. Since we are considering the periodic problem, $m$ must be
an even number and therefore, $m \geq 2(n+1).$  Also, note that
for any given even and natural number $q \geq 2(n+1),$ function
$b(x) \equiv \lambda_{q}$ belongs to $\landan$ and for function
$v(x) = \sin \frac{q\pi x}{T},$ we have $m = q.$ In this way, we
have proved the first two parts of the theorem.

Continuing with the proof,  if $i,$ with $ \ 0 \leq i \leq 2m-1$
is given and $u$ satisfies (\ref{2806076}), then
$$
\begin{array}{c}
\int_{t_i}^{t_{i+1}} u'^2(t) = \int_{t_i}^{t_{i+1}} a(t)u^2(t) =
\\ \\
\int_{t_i}^{t_{i+1}} (a(t)-\landaene) u^2 (t) +
\int_{t_i}^{t_{i+1}} \landaene u^2 (t)
\end{array}
$$
Therefore, $$ \int_{t_i}^{t_{i+1}} u'^2(t) - \landaene
\int_{t_i}^{t_{i+1}} u^2(t) \leq \Vert a - \landaene \Vert_{L^1
(t_i,t_{i+1})} \Vert u^2 \Vert_{L^\infty (t_i,t_{i+1})}
$$
Since $u'$ has no zeros in the interval $(t_i,t_{i+1})$ and
$u(t_i) = 0,$ we have $\Vert u^2 \Vert_{L^\infty(t_i,t_{i+1})} =
u^2 (t_{i+1}).$ This proves (see (\ref{050420111}))
\begin{equation}\label{202071}
\Vert a - \landaene \Vert_{L^1 (t_i,t_{i+1})} \geq
\frac{\int_{t_i}^{t_{i+1}} u'^2 - \landaene \int_{t_i}^{t_{i+1}}
u^2 }{u^2(t_{i+1})} \geq \frac{2n\pi}{T} \cot (\frac{2n\pi}{T}
(t_{i+1} -t_i)).
\end{equation}
In particular, this implies
\begin{equation}\label{2705091}
\Vert a - \landaene \Vert_{L^1 (0,T)} \geq \frac{2n\pi}{T}
\sum_{i=0}^{2m-1} \cot (\frac{2n\pi}{T} (t_{i+1} -t_i))
\end{equation}
The right-hand side of (\ref{2705091}) attains its minimum if and
only if $t_{i+1}-t_i = \frac{T}{2m}, \ 0 \leq i \leq 2m-1$ (see
Lemma 2.5 in \cite{cavijems}) so that
\begin{equation}\label{2805091} \gamma_{1,n} \geq T \lambda_{2n-1} + \frac{4n\pi}{T}m
\cot \frac{n\pi}{m}
\end{equation}
Taking into account that the function $ m \cot \frac{n\pi}{m}$ is
strictly increasing with respect to $ m,$ and that $m \geq
2(n+1),$ we deduce
\begin{equation}\label{03071}
\gamma_{1,n}  \ \geq T \lambda_{2n-1} + \frac{8\pi n(n+1)}{T} \cot
\frac{n\pi}{2(n+1)}.
\end{equation}

At this point, an appropriate minimizing sequence for
$\gamma_{1,n}$ may be constructed. The details are very technical
and they do not contribute anything important. On the other hand,
a careful observation of the optimality properties of the
different previos inequalities, do possible to prove that
$\gamma_{1,n}$ is not attained (see Lemma 2.3 and Lemma 2.4 in
\cite{cavijmaa}).

For the last part of the theorem, let us assume that the function
$a$ satisfies (\ref{positivo}). Then, since $\lambda_{2n-1} \prec
a,$ we trivially have $\lambda_{2n}(a) <
\lambda_{2n}(\lambda_{2n-1}) = 0.$ To prove that $\lambda_{2n+1}
(a) >0$ we use a continuation method: let us define the continuous
function $g: [0,1] \rightarrow \real$ by
$$
g(\varepsilon) = \lambda_{2n+1}(a_{\varepsilon}(\cdot))
$$
where $a_{\varepsilon}(\cdot) = \lambda_{2n-1} + \varepsilon
(a(\cdot) - \lambda_{2n-1}).$ Then $g(0) =\lambda_{2n+1}\left (
\lambda_{2n-1} \right )= \lambda_{2n+1} - \lambda_{2n-1} > 0.$
Moreover, $g(\varepsilon) \neq 0, \ \forall \ \varepsilon \in
(0,1].$ In fact, we deduce from the previous parts of the theorem
that the number $0$ is not an eigenvalue of the function
$a_{\varepsilon}$ for the periodic boundary conditions. As a
consequence, $\lambda_{2n+1}(a)= g(1)
> 0$ and the theorem is proved.

\end{proof}

\begin{remark}\label{010609t1} We can obtain similar results if, in the definition
of the set $\Lambda_n$ in (\ref{p2806073}), we consider $n \in
\real^+$ instead of $n \in \N$. This may be important in order to
obtain new stability results for periodic linear equations, as it
is shown in Theorem \ref{t0605091} below. Only some minor changes
are necessary. From this point of view, if we consider
$\gamma_{1,n}$ as a function of $n \in (0,+\infty),$ then $\lim_{n
\rightarrow 0^+} \ \gamma_{1,n} = \frac{16}{T},$ the constant of
the classical $L^1$ Lyapunov inequality for the periodic problem
at the first eigenvalue,  which was obtained in
\cite{huaizhongyongjde94} by using methods of optimal control
theory.

$ $

In the case $n = 0,$ we can use similar reasonings taking
\begin{equation}
{\small \Lambda_0 = \{ a \in \periodica \setminus \{ 0 \}: 0 \leq
\int_0^T a(t) \ dt \ \ \mbox{and} \ (\ref{p1}) \ \mbox{has
nontrivial solutions} \ \} }
\end{equation}
In this case $m\geq 2$ and any even value $m \geq 2$ is possible.
Consequently
\begin{equation}\label{1307094}
\gamma_{1,0} = \inf_{a \in \Lambda_0} \ \Vert a^+ \Vert_{L^1
(0,T)} = \displaystyle \frac{16}{T}
\end{equation}
Let us remark that the restriction
$$
a \in \periodica \setminus \{ 0 \}: \lambda_0 = 0 \leq \int_0^T
a(t) \ dt
$$
is more general that the restriction $ \lambda_0  \prec a$
\end{remark}
(As a consequence, it is possible to obtain a $L_1-L_1$ Lyapunov
inequality at the first eigenvalue).
\begin{remark}\label{1212073}
 The case where $T = 2\pi$ and
function $a$ satisfies the condition $ A \leq a(t) \leq B, \
\mbox{a.e. in} \ (0,2\pi) $ where $k^2 < A < (k+1)^2 < B$ for some
$k \in \natu \cup \{ 0 \}, $ has been considered in
\cite{wanglina95}, where the authors also use optimal control
theory methods.
\end{remark}

\begin{remark}\label{02032011}
We can use our methods to do an analogous study for other boundary
value problems (see the next theorem for antiperiodic boundary
conditions). Consider that, after seeing the above reasonings, the
key point is to have an optimal knowledge about the number and
distribution of zeros of the functions $u$ and $u',$ moreover of
knowing the best value of the constant $m$ (given in the previous
theorem). For example, for the periodic problem the optimal value
of $m$ is $2(n+1)$ while for the antiperiodic problem the optimal
value is $m = 2n+1.$ See, for example, \cite{cavijems} for the
case of Neumann, Dirichlet or mixed boundary conditions at higher
eigenvalues.
\end{remark}

$ $

Thinking in the next subsection, we finish this part with a
similar theorem to Theorem \ref{pt1} for the antiperiodic problem.
\begin{theorem}\label{antipt1}
Let $n \in \natu$ and $a \in \tilde{\Lambda}_{n}$ be given and $u$
any nontrivial solution of (\ref{ap1}) such that $u(0) = u(T) =
0.$ If the zeros of $u$ in $[0,T]$ are denoted by $0 = t_0 < t_2 <
\ldots < t_{2m} = T$ and the zeros of $u'$ in $(0,T)$ are denoted
by $ t_1 < t_3 < \ldots < t_{2m-1}, $ then:
\begin{enumerate}
\item $t_{i+1} - t_{i} \leq \frac{T}{2(2n-1)}, \ \forall \ i: \ 0
\leq i \leq 2m-1.$ Moreover, at least one of these inequalities is
strict. \item $m$ is an odd number and $m \geq 2n+1.$ Any odd
value $m\geq 2n+1$ is possible. %
\item $\tilde{\gamma}_{1,n}$ is not attained and
\begin{equation}\label{ap1009071}
\tilde{\gamma}_{1,n} \equiv \inf_{a \in \tilde{\Lambda}_{n} } \
\Vert a  \Vert_{L^1 (0,T)} = T\tilde{\lambda}_{2n-1} +
\frac{2\pi(2n-1)(2n+1)}{T} \cot \frac{(2n-1)\pi}{2(2n+1)}
\end{equation}
\item If $a \in \periodica$ satisfies
\begin{equation}\label{apositivo}
\tilde{\lambda}_{2n-1} \prec a, \ \Vert a \Vert_{L^1(0,T)} \leq
\tilde{\gamma}_{1,n},
\end{equation}
then
\begin{equation}\label{a1604096}
\tilde{\lambda}_{2n}(a) < 0 < \tilde{\lambda}_{2n+1}(a)
\end{equation}
\end{enumerate}
\end{theorem}

\begin{remark}\label{1307092} A similar theorem to the previous one may be proved if $a
\in \tilde{\Lambda}_0$ where
\begin{equation}
{\small \tilde{\Lambda}_0 = \{ a \in \periodica:  \ (\ref{ap1}) \
\mbox{has nontrivial solutions} \ \} }
\end{equation}
In this case $m\geq 1$ and any even value $m \geq 1$ is possible.
Consequently
\begin{equation}\label{1307093}
\tilde{\gamma}_{1,0} = \inf_{a \in \tilde{\Lambda}_0} \ \Vert a^+
\Vert_{L^1 (0,T)} = \displaystyle \frac{4}{T}
\end{equation}
Let us remark that the restriction $ 0 \prec a $ which is natural
for the periodic problem (\ref{p1}), is not necessary in this
case. This additional restriction is necessary in the Neumann
problem (see Remark 4 in \cite{camovimia}), but not in the case of
Dirichlet or mixed boundary conditions.
\end{remark}

\subsection{Stability of linear periodic equations}

In this part of the paper, we deal with some stability properties
of the Hill's equation
\begin{equation}\label{nuevahille}
u''(t) + a(t) u(t) = 0, \ a \in \periodica.
\end{equation}
As we said in the Introduction, it is convenient to introduce the
parametric equation
\begin{equation}\label{2hilleq}
u''(t) + (\mu+a(t)) u(t) = 0, \ a \in \periodica, \ \mu \in \real.
\end{equation}
If $\lambda_i (a), \ i \in \N \cup \{ 0 \}$ and $\tilde{\lambda}_i
(a), \ i \in \N,$ denote, respectively the eigenvalues of
(\ref{nuevahille}) for the periodic and antiperiodic problem, then
it is known (\cite{cole}, \cite{hale}) that
\begin{equation}\label{2004097}
\lambda_0 (a) < \tilde{\lambda}_1 (a) \leq \tilde{\lambda}_2 (a) <
\lambda_{1}(a) \leq \lambda_{2}(a) < \tilde{\lambda}_{3}(a) \leq
\tilde{\lambda}_{4}(a) < \lambda_3 (a) \leq \ldots
\end{equation}
and that equation (\ref{2hilleq}) is stable if
\begin{equation}\label{2004095}
\mu \in (\lambda_{2n}(a),\tilde{\lambda}_{2n+1}(a)) \cup
(\tilde{\lambda}_{2n+2}(a),\lambda_{2n+1}(a))
\end{equation}
for some $n \in \N \cup \{ 0 \}$ and that equation (\ref{2hilleq})
is unstable if
\begin{equation}\label{2004096}
\mu \in (-\infty,\lambda_0 (a)] \cup
(\lambda_{2n+1}(a),\lambda_{2n+2}(a)) \cup
(\tilde{\lambda}_{2n+1}(a),\tilde{\lambda}_{2n+2}(a))
\end{equation}
for some $n \in \N \cup \{ 0 \}.$ If $\mu = \lambda_{2n+1}(a)$ or
$\mu = \lambda_{2n+2}(a)$, (\ref{2hilleq}) is stable if and only
if $\lambda_{2n+1}(a) = \lambda_{2n+2}(a)$ and, finally, if $\mu =
\tilde{\lambda}_{2n+1}(a)$ or $\mu = \tilde{\lambda}_{2n+2}(a)$,
(\ref{2hilleq}) is stable if and only if
$\tilde{\lambda}_{2n+1}(a) = \tilde{\lambda}_{2n+2}(a).$

As a consequence, if we establish conditions which assure either
$\lambda_{2n}(a) < 0 < \tilde{\lambda}_{2n+1}(a)$ or
$\tilde{\lambda}_{2n+2}(a) < 0 < \lambda_{2n+1}(a)$ then
(\ref{nuevahille}) is stable. In this case, we say that $\mu = 0$
belongs to the $n^{th}$ stability zone of (\ref{nuevahille}).

$ $

In what follows in this subsection, we choose $T = \pi,$ for
simplicity. If one uses $L_\infty$ Lyapunov inequalities, the
following result may be proved ( see \cite{magnus}, Chapter V,
Theorem 5.5):

If $r$ and $s$ are given real numbers and
\begin{equation}\label{2905091}
r^2 \leq  a(t) \leq s^2
\end{equation}
then (\ref{nuevahille}) is stable for all possible functions $
a(\cdot)$ satisfying (\ref{2905091}) if and only if the interval
$(r^2,s^2)$ does not contain the square of an integer.

In particular, concerning to the first stability zone,
(\ref{nuevahille}) is stable if
\begin{equation}\label{2905092}
0 \leq a(t) \leq  1
\end{equation}
and for functions satisfying $0 \leq a(t),$ this result is optimal
in the following sense: for any positive number $\varepsilon$
there is some function $ a(t)$ with $a \in \periodica,$ satisfying
$0 \leq a(t) \leq 1+\varepsilon$ and such that (\ref{nuevahille})
is unstable.

$ $

The previous criterion on the stability of (\ref{nuevahille}) may
be considered as a $L_\infty-L_\infty$ criterion. We can exploit
the results obtained in Theorem \ref{pt1} and Theorem
\ref{antipt1} to obtain new results on the stability properties of
(\ref{nuevahille}) of the type $L_\infty-L_1$. This is the purpose
of the next theorem, where function $a$ can be cross an arbitrary
number of squares of an integer, as long as its $L_1$-norm is
properly controlled (see \cite{cavijmaa} for a more detailed
proof).

\begin{theorem}\label{t0605091}
Let $a \in \periodicapi$  satisfying
\begin{equation}\label{borg}
\begin{array}{c}
 \exists \ p \in \N, \ \exists \ k \ \in [p^2,(p+1)^2]: \\ \\
 k \leq  a, \ \
 \Vert a \Vert_{L^1(0,\pi)}  \leq k\pi + k^{1/2}2(p+1)\cot \frac{k^{1/2}\pi}{2(p+1)}
 \end{array}
\end{equation}
Then $\mu = 0$ is in the $n^{th}$ stability zone of the Hill's
equation (\ref{2hilleq}).
\end{theorem}

\begin{proof} It is clearly not restrictive (see Theorem \ref{pt1}) to assume that
\begin{equation}\label{bisborg}
\begin{array}{c}
 \exists \ p \in \N, \ \exists \ k \ \in (p^2,(p+1)^2): \\ \\
 k \prec  a, \ \
 \Vert a \Vert_{L^1(0,\pi)}  \leq k\pi + k^{1/2}2(p+1)\cot \frac{k^{1/2}\pi}{2(p+1)}
 \end{array}
\end{equation}
Let us suppose, for instance, that $p = 2n, \ n \in \N.$ Then,
$\lambda_{2n}(a) < \lambda_{2n}(\lambda_{2n-1}) = 0.$ On the other
hand, since $k \prec a,$ doing a similar reasoning to that in
Theorem \ref{pt1}, but for the antiperiodic problem, we have that
if $u$ is a nontrivial solution of (\ref{ap1}) such that $u(0) =
u(\pi) =0,$ then $\vert x_{i+1}-x_i \vert \leq
\frac{\pi}{2k^{1/2}}.$ This implies the relation $m
> k^{1/2}$ in in Theorem \ref{antipt1}. But since we are now considering the antiperiodic problem
(\ref{ap1}), $m$ must be an odd number. Also $p < k^{1/2} < p+1,$
and as $p =2n,$ we deduce $m \geq 2n+1.$ Consequently,
\begin{equation}\label{1105091}
\begin{array}{c}
\Vert a - k \Vert_{L^1 (0,\pi)} \geq  k^{1/2}2(2n+1)\cot
(\frac{k^{1/2}\pi}{2(2n+1)})
\end{array}
\end{equation}
(see (\ref{2805091})).  Moreover, this last constant is not
attained. As in Theorem \ref{pt1}, if  $h: [0,1] \rightarrow
\mathbb{R}$ is defined as por $h(\varepsilon) =
\tilde{\lambda}_{2n+1}\left ( k + \varepsilon (a(\cdot) - k)
\right ),$ we obtain $h(0) > 0$ and (from (\ref{bisborg})),
$h(\varepsilon) \neq 0, \ \forall \ \varepsilon \in (0,1].$ Then,
$h(1) = \tilde{\lambda}_{2n+1}(a) > 0.$ As a consequence, $\mu = 0
\in (\lambda_{2n}(a),\tilde{\lambda}_{2n+1}(a))$ and the theorem
is proved. The proof is similar if $p$ is an odd number.
\end{proof}

\begin{remark} The case where $a(t) = \alpha + \beta \psi (t),$
with $\psi \in \periodicapi, \ \int_0^\pi \psi (t) \ dt = 0$ and
$\int_0 ^\pi \vert \psi (t) \vert \ dt = 1/\pi,$ was studied by
Borg (\cite{borg}). Borg used the characteristic multipliers
determined from Floquet's theory. He deduced stability criteria
for (\ref{2hilleq}) by using the two parameters $\alpha$ and
$\beta.$ For a concrete function $a,$ this implies the use of the
two quantities
$$
\int_0^\pi a(t) \ dt, \ \  \left \Vert a(\cdot) - \frac{1}{\pi}
\int_0^\pi a(t) \ dt \right \Vert_{L^1 (0,\pi)}
$$
It is clear that the results given in Theorem \ref{t0605091} are
of a different nature (see \cite{magnus} and the translator's note
in \cite{krein}). In fact, our results are similar to those
obtained by Krein \cite{krein} by using a different procedure.
However, Krein assumed $k = p^2$ and an strict inequality for
$\Vert a \Vert_{L^1(0,\pi)}$ in (\ref{borg}) (see Theorem 9 in
\cite{krein}). By using Theorem \ref{antipt1} we can assume a non
strict inequality in (\ref{borg}) since the constant
$\tilde{\gamma}_{1,n}$ is not attained.

Finally, if for a given function $a \in \periodicapi$ we know that
$a$ satisfies (\ref{bisborg}), the result given in Theorem
\ref{t0605091} is more precise than Krein's result since the
function
$$
k\pi + k^{1/2}2(p+1)\cot \frac{k^{1/2}\pi}{2(p+1)}, \ k \in \
[p^2,(p+1)^2]
$$
is strictly increasing with respect to $k.$
\end{remark}

\section{Systems of ordinary differential equations}

In this section we present some results on Lyapunov constants and
stability properties for systems of equations
\begin{equation}\label{2ieq35}
u''(t) + P(t)u(t)= 0,\ t \in \real,
\end{equation}
where the matrix function $P(\cdot)$ is $T-$periodic.

In regard to stability properties of (\ref{2ieq35}), the results
proved by Krein in \cite{krein} show that the problem is closely
related to Lyapunov inequalities. In \cite{krein}, the author
assumes that $P(\cdot) \in \Lambda,$ where $\Lambda$ is defined as

\afirm{[$\Lambda$]}{10cm}{The set of real $n\times n$ symmetric
matrix valued function $P(\cdot)$, with continuous and
$T-$periodic element functions $p_{ij}(t), \ 1 \leq i,j \leq n,$
such that (\ref{2ieq35}) has not nontrivial constant solutions and
$$ \intl \langle P(t)k,k \rangle  \ dt \geq 0, \
\forall \ k \in \real^n. $$ }

Here, $\langle \cdot , \cdot \rangle$ denotes the usual inner
product in $\R^n.$

Let us remark that, in the scalar case, the previous hypotheses
imply that $\lambda_0 (P) < 0,$ where $\lambda_0 (P)$ is the first
eigenvalue for the periodic problem. Therefore, and as in the
scalar case, to have some stability property for the system
({\ref{2ieq35}), it is necessary another condition involving the
antiperiodic boundary value problem.

$ $

The system (\ref{2ieq35}) is said to be stably bounded
(\cite{krein}) if there exists $\varepsilon = \varepsilon (P) \in
\real ^+,$ such that all solutions of the system
\begin{equation}\label{2eq35}
u''(t) + Q(t)u(t)= 0, \ t \in \real,
\end{equation}
are bounded for all matrix function $Q(\cdot) \in \Lambda,$
satisfying
$$
\displaystyle \max_{1 \leq i,j \leq n} \ \intl \vert p_{ij}(t) -
q_{ij}(t) \vert \ dt < \varepsilon,
$$
that is, for $L_1$ small perturbations of the matrix function
$P(\cdot).$

Krein proved that if $P(\cdot) \in \Lambda,$ then the system
(\ref{2ieq35}) is stably bounded if $\lambda_1
> 1,$ where $\lambda_1$ is the smallest positive eigenvalue of the
antiperiodic eigenvalue problem
\begin{equation}\label{2eq37}
u''(t) + \lambda P(t)u(t)= 0, \ t \in \real, \ u(0)+u(T) =
u'(0)+u'(T) = 0.
\end{equation}
The result given by Krein is a nice theoretical result, but for
systems of equations, and assuming $P(\cdot) \in \Lambda,$ it is
not easy to give sufficient conditions to ensure the property
$\lambda_1 > 1.$ In this section, we establish some new conditions
to get this. Previously, we need to characterize the best $L_p$
Lyapunov constants as the minimum of some special minimization
problems (constrained or unconstrained, depending on the boundary
value problem that we are considering).

A key point is that these minimization problems involve only
appropriate functionals defined on some subsets of the
corresponding Sobolev space, without any reference to the function
$a(\cdot)$ considered in the definition of the Lyapunov constants
(for instance, for the case $p=1$, see the definition of the
constants $\gamma_{1,n}$ and $\tilde{\gamma}_{1,n}$ in
(\ref{pconstante}) and (\ref{anticonstante}), respectively).

\subsection{$L_p$ Lyapunov constants and minimization problems for ODE}

This subsection will be concerned with some preliminary results on
Lyapunov inequalities for the antiperiodic boundary value problem
\begin{equation}\label{eq2}
u''(t) + a(t)u(t) = 0, \ t \in (0,T), \ u(0)+u(T) = u'(0) + u'(T)
= 0,
\end{equation}
where $a \in L_{T} (\R,\R).$ We consider the antiperiodic problem
because we will use these results in the next subsection, but the
same idea may be applied to other boundary value problems with
Neumann, mixed, Dirichlet, periodic conditions, etc. (see
\cite{camovimia}, \cite{camovim} \cite{cavijmi}, \cite{cavina},
\cite{cavijde})).

As in the case $p=1,$ if we define the sets

\begin{equation}\label{eq4}
\Lambda^{ant} = \{ a \in \periodica : \ (\ref{eq2}) \ \mbox{has
nontrivial solutions} \ \},
\end{equation}
then for each $p$ with $1 \leq p \leq \infty,$ we can define the
$L^p$ Lyapunov constant $\beta_p ^{ant}$ for the antiperiodic
problem, as the real number
\begin{equation}\label{eq5}
\beta_{p}^{ant} \equiv \displaystyle \inf_{a \in
\Lambda^{ant}\bigcap L^p (0,T)} \ \ \Vert a^+ \Vert_p
\end{equation}
where
\begin{equation}\label{eq6}
\begin{array}{c}
\Vert a^+ \Vert_{p} = \left ( \displaystyle \int_{0}^{T} \vert a^+
(t)\vert ^{p} \ dt \right ) ^{1/p}, 1 \leq p < \infty; \ \ \Vert
a^+ \Vert_{\infty} =  sup \ ess \ \ a^+.
\end{array}
\end{equation}
An explicit expression for the constant $\betaant $, as a function
of $p$ and $T,$ has been obtained in \cite{zhang} (see also
\cite{camovimia}, \cite{cavidcds} and \cite{talenti} for the case
of Neumann, mixed and Dirichlet boundary conditions,
respectively). As mentioned earlier, in the next theorem we
characterize the $L^p$ Lyapunov constant $\beta_{p}^{ant},$ as a
minimum of a convenient minimization scalar problem. This will
allow the extension of the results to systems of equations. On the
other hand, the next theorem is fundamental to intuit the kind of
the results we can expect in PDE which will be considered in the
next section.

\begin{theorem}\label{l2}
If $1 \leq p \leq \infty$ is a given number, let us define the
sets $X_p ^{ant}$ and the functional $I_p ^{ant}: X_p
^{ant}\setminus \{ 0 \} \rightarrow \real,$ as
\begin{equation}\label{eqant10}
\begin{array}{c}
X_p ^{ant} = \left \{ v \in \sobol: v(0)+ v(T)= 0 \right \}, \ 1
\leq p \leq \infty, \\ \\
I_1 ^{ant}(v) = \displaystyle \frac{\dis \intl v'^{2}}{\Vert v
\Vert_{\infty}^{2}}, \ \ I_{\infty} ^{ant} (v) = \displaystyle
\frac{\dis \intl v'^{2}}{\displaystyle \intl v^{2}} \\ \\
I_p ^{ant}(v) = \dis \frac{\dis \intl v'^{2}}{\left ( \dis \intl
\vert v \vert ^{\frac{2p}{p-1}} \right )^{\frac{p-1}{p}}}, \
\mbox{if} \ \ 1 < p < \infty.
\end{array}
\end{equation}
Then, the $L_p$ Lyapunov constant $\beta_p ^{ant}$ defined in
(\ref{eq5}), satisfies
\begin{equation}\label{eqant11} \begin{array}{c}
\beta_{p}^{ant} = \displaystyle \min_{X_p ^{ant} \setminus \{ 0
\}} \ I_{p} ^{ant},\ \ 1 \leq p \leq \infty.
\end{array}
\end{equation}
\end{theorem}
\begin{proof}
In the proof, only the case $1 < p < \infty$ is considered (see
\cite{cavina} for the other cases). In summary, we apply the
following ideas:
\begin{enumerate}
\item The demonstration that the infimum of the functional $I_{p}
^{ant}$ is attained on $X_p ^{ant} \setminus \{ 0 \}$ is standard:
we prove the existence of bounded minimizing sequences and then we
use the weak lower semicontinuity property of the functional $I_p
^{ant}$. \item Once we have proved the existence of global
minimum, we apply the Lagrange multiplier theorem to obtain the
second order (nonlinear) differential equation which must be
satisfied by the function $u_0 \in X_p ^{ant} \setminus \{ 0 \}$
where the minimum of $I_{p} ^{ant}$ is attained (Euler equation).
\item Finally, we integrate the mentioned Euler equation. This is
not an easy task and here, again, the number and distribution of
the zeros of $u_0$ and $u_0 '$ is of an important help.
\end{enumerate}

To carry out the above ideas, let us denote
$$M_{p} = \inf_{X_p ^{ant}\setminus \{ 0 \}}\ I_p ^{ant}.$$ If $\{ u_{n} \} \subset
X_{p}^{ant} \setminus \{ 0 \}$ is a minimizing sequence, since the
sequence $\{ k_{n}u_{n} \}, \ k_{n} \neq 0, $ is also a minimizing
sequence, we can assume without loos of generality that $\dis
\intl \vert u_{n} \vert ^{\frac{2p}{p-1}} = 1.$ Then $\left \{
\dis \intl \vert u_{n}'^{2} \vert \right \} $ is also bounded.
Moreover, for each $u_{n}$ there is $x_{n} \in [0,T]$ such that
$u_{n}(x_{n}) = 0.$ Therefore, $\{ u_{n} \}$ is bounded in
$\sobol.$ So, we can suppose, up to a subsequence, that $u_{n}
\rightharpoonup u_{0}$ in $\sobol$ and $u_{n} \rightarrow u_{0}$
in $C[0,L]$ (the space of continuous functions in $[0,L]$ with the
uniform norm). The strong convergence in $C[0,L]$ gives us $u_0
(0)+u_0 (T) = 0$. Therefore, $u_{0} \in X_{p}^{ant} \setminus \{ 0
\}.$ The weak convergence in $\sobol$ implies $I^{ant}_{p}(u_{0})
\leq \liminf \ I^{ant}_{p}(u_{n}) = M_{p}.$ Then $u_{0}$ is a
minimizer.

$ $

We conclude,
$$
H'(u_{0})(v) = 0, \ \forall v \in \sobol \ \mbox{such that} \
v(0)+v(T) = 0.
$$
Here $H: \sobol \rightarrow \real$ is defined by
$$
H(u) = \dis \intl u'^{2} - M_{p} \left ( \dis \intl \vert u \vert
^{\frac{2p}{p-1}} \right )^{\frac{p-1}{p}}
$$
This implies that $u_0$ satisfies the Euler equation
\begin{equation}\label{eqant13}
\begin{array}{c}
u_0''(t) + A_{p}(u_0)\vert u_0 (t)\vert^{\frac{2}{p-1}}u_{0}(t) =
0, \ t \in (0,T), \\  u_0 (0)+u_0 (T) = 0, \ u_0 '(0)+u_0 ' (T) =
0,
\end{array}
\end{equation}
where
\begin{equation}\label{eqant14}
A_{p}(u_0) = M_{p} \left ( \dis \intl \vert u_0\vert
^{\frac{2p}{p-1}} \right )^{\frac{-1}{p}}
\end{equation}
Since the function $a(t) \equiv A_{p}(u_0)\vert u_0
(t)\vert^{\frac{2}{p-1}}$ satisfies $a(0) = a(T),$ it is not
restrictive to assume that, additionally, $u_0 (0) = u_0 (T) = 0.$
Euler equation (\ref{eqant13}) can be integrated, using similar
ideas to those considered in \cite{camovimia} Lemma 2.7  and
\cite{camovim} Theorem 2.1, for the Neumann problem. Then, we
deduce
\begin{equation}\label{eqant17}
M_{p} = \dis \frac{4n^{2}I^{2}p}{T^{2-\frac{1}{p}}
(p-1)^{1-\frac{1}{p}} (2p-1)^{1/p}}, \end{equation} where $ \ I =
\dis \int_{0}^{1} \frac{ds}{\left ( 1 - s^{\frac{2p}{p-1}} \right
) ^{1/2}}$ and $n \in \N$ is such that we denote the zeros of
$u_0$ in $[0,T]$ by $0 = x_0 < x_2 < \ldots < x_{2n} = T$ and the
zeros of $u_0 '$ in $(0,T)$ by $ x_1 < x_3 < \ldots < x_{2n-1}.$

For the antiperiodic problem, (and for the Neumann and Dirichlet
problem), $n \geq 1$ and any value $n \geq 1$ is possible. The
conclusion is that
\begin{equation}\label{eqant18}
M_{p} = \dis \frac{4I^{2}p}{T^{2-\frac{1}{p}}
(p-1)^{1-\frac{1}{p}} (2p-1)^{1/p}},
\end{equation}
which is the same constant as in the Neumann and Dirichlet
problem. Finally, in \cite{zhang} it is shown that this is,
exactly, the $L_p$ Lyapunov constant for the antiperiodic problem.
The theorem is proved.
\end{proof}

\subsection{Stability for linear periodic systems}

In this subsection we present some $L_p$ conditions  ($1 \leq p
\leq \infty$) which allow to prove that the smallest positive
eigenvalue of the antiperiodic eigenvalue problem (\ref{2eq37})
satisfies $\lambda_1> 1.$  As a consequence, if $P(\cdot) \in
\Lambda,$ the system (\ref{2ieq35}) is stably bounded. These
conditions are given in terms of the $L^p$ norm of appropriate
functions $b_{ii}(t), \ 1 \leq i \leq n,$ related to
(\ref{2ieq35}) and (\ref{2eq37}) through the inequality $P(t) \leq
B(t), \ \forall\ t \in [0,T],$ where $B(t)$ is a diagonal matrix
with entries given by $b_{ii}(t), \ 1 \leq i \leq n.$ These
sufficient conditions are optimal in the sense explained in Remark
\ref{remark3} below. Here, the relation $C \leq D$ between $n
\times n$ symmetric matrices means that $D-C$ is positive
semi-definite.

$ $

From now on we assume that the matrix function $P(\cdot) \in
\Lambda$ ($\Lambda$ was defined at the beginning of this section).
In this case, $\lambda_1$ has a variational characterization (see
\cite{krein}) given by
\begin{equation}\label{eq38}
\displaystyle \frac{1}{\lambda_1} = \max_{y \in G_T}  \ \intl
\langle P(t)y(t),y(t)\rangle \  dt,
\end{equation}
where
\begin{equation}\label{39}
G_T = \{ y \in (H^1 (0,T))^n: y(0)+y(T) = 0, \ \dis \sum_{i=1}^{n}
\intl (y_i ' (t))^2 \ dt = 1 \}.
\end{equation}

\begin{theorem}\label{t2}
Let $P(\cdot) \in \Lambda$ be such that there exist a diagonal
matrix $B(t)$ with continuous and $T-$periodic  entries
$b_{ii}(t),$ and $p_i \in [1,\infty], \ 1 \leq i \leq n,$
satisfying
\begin{equation}\label{eq40}
\begin{array}{c}
P(t) \leq B(t), \ \forall \ t \in \real, \\ \\ \Vert b_{ii}^+
\Vert_{p_i} < \beta_{p_i}^{ant}, \ \mbox{if} \ p_i \in (1,\infty],
\ \ \Vert b_{ii}^+ \Vert_{p_i} \leq  \beta_{p_i}^{ant}, \
\mbox{if} \ p_i =1.
\end{array}
\end{equation}
Then, the system (\ref{2ieq35}) is stably bounded.
\end{theorem}
\begin{proof}
Let $y \in G_T.$ Then by using the Theorem \ref{l2} and more
specifically the relation (\ref{eqant11}), we have
\begin{equation}\label{eq41}
\begin{array}{c}
\intl \langle P(t)y(t),y(t)\rangle  \ dt \leq \intl \langle B(t)y(t),y(t)\rangle  \ dt \leq \\
\\
\dis \sum_{i=1}^{n} \ \dis \intl b_{ii}(t)(y_i (t))^2 (t) \ dt
\leq \dis \sum_{i=1}^{n} \Vert b_{ii}^+ (t) \Vert_{p_i} \Vert y_i
^2 \Vert_{\frac{p_i}{p_i -1}}\leq  \\ \\
\leq  \dis \sum_{i=1}^{n} \beta_{p_i}^{ant} \Vert y_i ^2
\Vert_{\frac{p_i}{p_i -1}} \leq \dis \sum_{i=1}^{n} \intl (y_i '
(t))^2 \ dt = 1, \ \forall \ y \in G_T,
\end{array}
\end{equation}
where
$$
\begin{array}{c}
\frac{p_i}{p_i -1} = \infty, \ \ \mbox{if} \ p_i = 1 \\ \\
\frac{p_i}{p_i -1} = 1, \ \ \mbox{if} \  p_i = \infty. \end{array}
$$
(\ref{eq41}) implies $\lambda_1 \geq 1.$ But if $\lambda_1 = 1,$
let us choose $y(\cdot)$ as any nontrivial solution of the problem
\begin{equation}\label{3eq37}
y''(t) + P(t)y(t)= 0, \ t \in \real, \ y(0)+y(T) = y'(0)+y'(T) =
0.
\end{equation}
Then some component, say $y_j,$ is nontrivial. If $p_j \in
(1,\infty],$ $(\beta_{p_j}^{ant}-\Vert b_{jj}^+ \Vert_{p_j})\Vert
y_j^2 \Vert_{\frac{p_j}{p_j -1}} >0 $ and
$(\beta_{p_i}^{ant}-\Vert b_{ii}^+ \Vert_{p_i})\Vert y_i^2
\Vert_{\frac{p_i}{p_{i} -1}}\geq 0, \ \forall \ i \neq j$, so that
we have a strict inequality in (\ref{eq41}). This is a
contradiction with (\ref{3eq37}). If $p_j =1,$ we can use a
similar reasoning (see \cite{cavina}). Consequently $\lambda_1 > 1
$ and the theorem is proved.
\end{proof}

\begin{remark}\label{remark3}
Previous theorem is optimal in the following sense: for any given
positive numbers $\gamma_i, \ 1 \leq i \leq n,$ such that at least
one of them, say $\gamma_j,$ satisfies
\begin{equation}\label{22optimalidad}
\gamma_j > \beta_{p_j} ^{ant}, \ \mbox{for some} \ p_j \in
[1,\infty],
\end{equation} there exists a diagonal  $n\times n$ matrix
$P(\cdot) \in \Lambda$ with entries $p_{ii}(t), \ 1 \leq i \leq
n,$ satisfying $\Vert p_{ii}^+ \Vert_{p_i} < \gamma_{i}, \ 1 \leq
i \leq n$ and such that the system (\ref{2ieq35}) is not stable.

$ $

To see this, if $\gamma_j$ satisfies (\ref{22optimalidad}), then
there exists some continuous and $T-$periodic function $p(t), $
not identically zero, with $\intl p (t) \ dt \geq 0, $ and $\Vert
p^+ \Vert_{p_j} < \gamma_j,$ such that the scalar problem
$$w''(t) + p(t)w(t) = 0,
$$
is not stable (see Theorem 1 in \cite{zhangli}). If we choose
$$
p_{jj}(t) = p(t), \ \ p_{ii}(t) = \delta \in \real^+, \ \mbox{if}
\  i\neq j,
$$
with $\delta$ sufficiently small, then (\ref{2ieq35}) is unstable.
\end{remark}

\begin{example}\label{ejemplo240510t} In the two
dimensional case, we show an example where the hypotheses of the
Theorem \ref{t2} may be checked directly by using the elements
$p_{ij}$ of the matrix $P(t).$ The example is based on a similar
one shown by the authors in \cite{cavijde}, in the study of
Lyapunov inequalities for elliptic systems (see also Example
\ref{c2p} below).

$ $

Let the matrix $P(t)$ be given by
\begin{equation}\label{teq10}
P(t) = \left (
\begin{array}{cc} p_{11}(t) & p_{12}(t) \\
p_{12}(t) & p_{22}(t) \end{array} \right )
\end{equation}
where \afirm{[{\bf H1}]}{10cm} {
$$
\begin{array}{c} p_{ij} \in C_T(\real,\real), \ 1 \leq i,j \leq 2, \\ \\  p_{11}(t) \geq 0, \
p_{22}(t) \geq 0, \ \det \ P(t) \geq 0, \ \forall \ t \in \real,
\\ \\ \det \ P(t) \neq 0, \ \mbox{for some}\ t \in \real.
\end{array}
$$ }
Here $C_T (\real,\real)$ denotes the set of real, continuous and
$T-$periodic functions defined in $\real$ and $\det P(t)$ means
the determinant of the matrix $P(t)$.

In addition, let us assume that there exist $p_1, p_2 \in
(1,\infty]$ such that
\begin{equation}\label{teq11}
\Vert p_{11}\Vert_{p_1} < \beta_{p_1}^{ant}, \ \ \Vert p_{22} +
\displaystyle \frac{p_{12}^2}{\beta_{p_1}^{ant} - \Vert p_{11}
\Vert_{p_1}} \Vert _{p_2} < \beta_{p_2}^{ant}.
\end{equation}
Then (\ref{2ieq35}) is stably bounded (see \cite{cavina}) for the
details).

$ $

Let us observe that from (\ref{teq11}) we deduce
\begin{equation}\label{2teq11}
\Vert p_{11}\Vert_{p_1} < \beta_{p_1}^{ant}, \ \ \Vert p_{22}
\Vert _{p_2} < \beta_{p_2}^{ant}.
\end{equation}
As a consequence, the uncoupled system
\begin{equation}\label{14092010}
v''(t) + R(t)v(t) = 0, \ t \in \real,
\end{equation}
where
\begin{equation}\label{2eq10}
R(t) = \left (
\begin{array}{cc} p_{11}(t) & 0 \\
0 & p_{22}(t) \end{array} \right )
\end{equation}
is stably bounded. Therefore, by using the definition of stably
bounded system, (\ref{2ieq35}) is stably bounded for any
continuous and $T-$periodic function $p_{12}$ with sufficiently
small $L_{1}-$ norm. However, (\ref{teq11}) does not imply,
necessarily, that the $L_{1}-$norm of the function $p_{12}$ is
necessarily small.
\end{example}

\section{Scalar partial differential equations}

In this section we consider the linear problem with Neumann
boundary conditions
\begin{equation}\label{ceq6}
\left.
\begin{array}{cl}
\Delta u(x)+ a(x)u(x)= 0, & x\in \Omega \, \\
 \frac{\partial u}{\partial n}(x)=0, \ \ \ \ \ \ \ &  x\in \partial \Omega
\end{array}
\right\}
\end{equation}
\noindent where $\Omega \subset \real^N$ ($N\geq 2$) is a bounded
and regular domain and  $\dis \frac{\partial}{\partial n}$ is the
outer normal derivative on $\partial \Omega.$ Our methods are also
applicable to the case of Dirichlet boundary conditions in
(\ref{ceq6}) and similar results can be obtained.

However, with respect to ODE there is an important difference: in
ODE, once we have solved a concrete problem, for example for
Dirichlet boundary conditions, by using appropriate changes of
variables we can obtain the best Lyapunov constants for other
boundary conditions such as Neumann, periodic, antiperiodic, mixed
boundary value problems, etc. (see \cite{zhang}). To the best of
our knowledge, we do not know similar results  for PDE and
consequently, the problem must be solved differently. We overcome
this difficulty characterizing the best constant as the minimum of
some convenient minimization problems as it was done for ODE in
Theorem \ref{l2} for the case of antiperiodic boundary conditions.

For Neumann problems like (\ref{ceq6}) the characterization must
be established by using some especial constrained minimization
problems, a fact that does not occur in the case of Dirichlet
boundary conditions (see Lemma 4.2 and Remark 5 in
\cite{camovijfa} and Theorem \ref{miap} below).

First we treat the scalar linear problem (\ref{ceq6}). Then we
apply the linear study to nonlinear resonant problems by using the
Schauder fixed point theorem (see Theorem \ref{nolineal}).
Analogous results can be obtained for nonlinear resonant problems
for ODE (see \cite{camovimia}, \cite{cavidcds}, \cite{cavijems}
and  \cite{cavijmaa}).

\subsection{Lyapunov constants for scalar partial differential equations}

\noindent This subsection will be concerned with the linear
boundary value problem (\ref{ceq6}). Remember that $\Omega \subset
\real^N$ is a bounded and regular domain, and we assume $a \in
\Sigma$, where
\begin{equation}\label{0803111}\Sigma = \{ \   a \in L^{\infty}(\Omega)\setminus \{0\}:
\dis \int_\Omega a(x) \ dx \geq 0 \ \mbox{and} \ (\ref{ceq6}) \
\mbox{has nontrivial solutions}\}
\end{equation}
As in the ordinary case, the positive eigenvalues of the
eigenvalue problem
\begin{equation}\label{ceq11}
\left.
\begin{array}{cl}
\Delta u(x)+\lambda u(x)= 0, & x\in \Omega \, \\
 \dis \frac{\partial u(x)}{\partial n}=0,\ \ \ \ \ \ \ &  x\in \partial \Omega
\end{array}
\right\}
\end{equation}
belong to $\Sigma$ and therefore, the quantity
\begin{equation}\label{1603113}
 \beta_{p} \equiv \inf_{a \in \Sigma } \ \Vert
a^+\Vert_p \ , \ 1 \leq p \leq \infty
\end{equation} is well
defined.

\begin{theorem}\label{t2partial} The following statements hold:
\begin{enumerate}
\item $\beta_{1} = 0$  and $\beta_{\infty} = \lambda_{1}, \
\forall \ N \geq 2.$ Here $\lambda_{1}$ denotes the first positive
eigenvalue of the eigenvalue problem (\ref{ceq11}). \item If
$N=2,$ $\beta_p>0 \Leftrightarrow p \in (1,\infty]$.

\noindent If $N\geq 3,$ $\beta_p>0 \Leftrightarrow p \in
[\frac{N}{2},\infty].$

\noindent If $N\geq 2$ and $\frac{N}{2}<p\leq \infty,$ then
$\beta_p$ is attained.

\item The mapping $(\frac{N}{2},\infty) \rightarrow \real,$ $p
\mapsto \beta_{p},$ is continuous and the mapping
$[\frac{N}{2},\infty] \rightarrow \real, $ $p \mapsto \vert \Omega
\vert^{-1/p}\beta_{p},$ is strictly increasing. Here, $\vert
\Omega \vert$ denotes the Lebesgue measure of the domain $\Omega.$

\end{enumerate}
\end{theorem}

\begin{proof}  {\it We summarize the main ideas}:
\begin{enumerate}
\item  If $N \geq 3$ and $\frac{N}{2} < p < \infty,$ the ideas are
the same as in the ordinary case (see Theorem 2.1 in
\cite{camovimia}). In fact, since $\frac{N}{2} < p,$ then $\dis
\frac{2p}{p-1} < \dis \frac{2N}{N-2}$ and consequently, the
embedding of the Sobolev space
$H^{1}(\Omega)$ into $L^{2p/p-1}(\Omega)$ is compact. %

 \item If $N \geq 3$ and $p = \frac{N}{2},$ then $
\frac{2p}{p-1} = \frac{2N}{N-2}$ and the embedding $H^{1}(\Omega)
\subset L^{2N/N-2}(\Omega)$ is continuous but not compact. This
implies that the infimum $\beta_{p}
> 0$, but we do not know if $\beta_{p}$ is a minimum.

\item If $N \geq 3$ and $1\leq p<\frac{N}{2}$, we prove that
$\beta_p=0$ by finding appropriate minimizing sequences. Roughly
speaking, the main idea is to take first a function $u$ and to
calculate the corresponding function $a$ for which $u$ is a
solution of $$ \Delta u(x)+ a(x)u(x) = 0, x\in \Omega,  \ \ \
\frac{\partial u(x)}{\partial n}=0, \ \  x\in
\partial \Omega
$$ Obviously, if $u$ is smooth enough, then we
must impose two conditions: i) $\frac{\partial u}{\partial n}=0$
on $\partial \Omega$, ii) The zeros of $u$ are also zeros of
$\Delta u$. For instance, if $\Omega=B(0,1)$ we can take radial
functions $u(x)=f(\vert x\vert)$ of the form $f(r)=\alpha
r^{-\delta}-\beta r^{-\gamma}$, ($\delta >0, \ \gamma > 0, 0 < r <
1$).

$ $

The case of general bounded regular domains $\Omega$ can be
treated as follows: first of all, note that if we define
$\Omega+x_0=\left\{ x+x_0\, ;\, x\in \Omega \right\}$ (for
arbitrary $x_0\in \real^N$), then
$\beta_p\left(\Omega+x_0\right)=\beta_p\left(\Omega\right)$. On
the other hand, if we define $r \, \Omega=\left\{r x\, ;\, x\in
\Omega \right\}$ (for arbitrary $r\in \real^+$), then
$\beta_p\left(r \, \Omega\right)=
r^{\frac{N}{p}-2}\beta_p\left(\Omega\right)$. Hence
$\beta_p(\Omega)=0\Leftrightarrow \beta_p\left(r \,
\Omega+x_0\right)=0$. Then, we can suppose without loss of
generality that $\overline{B}(0,1)\subset \Omega$.%

\item If $N=2$ and $p = 1,$  the ideas are the same as in the
previous step, but we use the fundamental solution $\ln \vert x
\vert$ to find appropriate minimizing sequences (see Lemma 3.2 in
\cite{camovijfa}).

\item If $N = 2,$ the embedding $H^{1}(\Omega) \subset
L^{q}(\Omega)$ is compact $\forall \ q \in [1,\infty)$ and
therefore, if $1 < p < \infty,$
the ideas are the same as in the ordinary case. %

\item The case $p=\infty$ is trivial if we use the variational
characterization of $\lambda_1.$

 \end{enumerate}

\end{proof}

\begin{remark}
The set $\Sigma$ defined in (\ref{0803111}) can be replaced by

If $N \geq 3,$
$$ \Sigma = \{ a \in
L^{\frac{N}{2}}(\Omega)\setminus \{0\}: \dis \int_\Omega a(x) \ dx
\geq 0 \ \mbox{and} \ (\ref{ceq6}) \ \mbox{has nontrivial
solutions}\} $$ and if $N=2,$
$$\Sigma = \{ a:\Omega \rightarrow \real \mbox{ s. t. } \exists
q\in(1,\infty]\\ \ \mbox{with} \ a\in L^q(\Omega)\setminus \{0\},
\
 \dis \int_\Omega a(x) \ dx \geq 0 \ \mbox{and} \ (\ref{ceq6})$$

\noindent\ has nontrivial solutions$\}$

\end{remark}

\begin{remark}\label{080311}

In the definition of the set $\Sigma$ we have imposed $\int_\Omega
a\geq 0$. This is not a technical but a natural assumption for
Neumann boundary conditions (an also for periodic boundary
conditions in the ordinary case. See Section 2 in \cite{cavina}).
Otherwise, the corresponding infimum will be always zero. To see
this, note that if $u\in H^1(\Omega)$ is a positive nonconstant
solution of (\ref{ceq6}) and we consider $v=\displaystyle
\frac{1}{u}$ as test function in the weak formulation, we obtain

$$
\displaystyle \int_\Omega \nabla u \cdot \nabla \left(
\frac{1}{u}\right) = \displaystyle \int_\Omega a\, u \frac{1}{u},
$$

\noindent which implies

$$
\displaystyle \int_\Omega a= \displaystyle - \int_\Omega
\frac{\vert \nabla u\vert^2 }{u^2}< 0.
$$

With this in mind, if we take a nonconstant $u_0\in C^2
(\overline{\Omega})$ such that $\frac{\partial u_0}{\partial
n}(x)=0, \, \forall x\in \partial \Omega$ then, for large $n\in
\natu$, we have that $u_n=u_0+n$ is a positive nonconstant
solution of (\ref{ceq6}), with $a_n=\frac{-\Delta u_0}{u_0+n}$.
Clearly $\displaystyle \int_\Omega a_n< 0$ and $\Vert a_n
\Vert_p\rightarrow 0$ as $n\rightarrow \infty$ for every $1\leq
p\leq \infty$.

\end{remark}

\begin{remark}
A property similar to the last part of the previous theorem is
valid for ODE (see, for instance, Theorem 2.1 in \cite{camovimia}
for the case of Neumann boundary conditions).
\end{remark}

\subsection{Resonant nonlinear scalar problems for PDE}

The results that we have obtained for linear problem (\ref{ceq6})
can be useful for the study of nonlinear boundary value problems
of the form
\begin{equation}\label{ceq9}
\left.
\begin{array}{cl}
\Delta u(x) + f(x,u(x))= 0, & x\in \Omega \, \\
 \frac{\partial u}{\partial n}(x)=0, \ \ \ \ \ \ \ &  x\in \partial \Omega
\end{array}
\right\}
\end{equation}
In fact, this is common in analysis: to study nonlinear problems
with the help of convenient linearized problems (implicit function
theorem, inverse function theorem, bifurcation theory, etc.)

$ $

\noindent As in the previous subsection we assume that $\Omega
\subset \real^N$ ($N\geq 2$) is a bounded and regular domain and
the function $f: \comega \times \real \rightarrow \real, \ (x,u)
\mapsto f(x,u),$ satisfies the condition

\medskip {\bf (H)} $f, f_{u}$ are Caratheodory functions
i.e. $f(x, u), f_{u}(x,u)$ are measurable in $x$ for every $u$ and
continuous in $u$ for almost every $x$ and $0 \leq f_{u}(x,u)$ in
$\comega\times \real.$

\medskip The existence of a solution of (\ref{ceq9}) implies
\begin{equation} \label{ceq10}
\dis \int_{\Omega} f(x,s_0) \ dx = 0
\end{equation}
for some $s_0 \in \real$. Trivially, conditions (H) and
(\ref{ceq10}) are not sufficient for the existence of solutions of
(\ref{ceq9}). Indeed, consider the problem
\begin{equation}\label{180305}
\left.
\begin{array}{cl}
-\Delta u(x)=\lambda_{1}u(x)+ \varphi_{1}(x), & x\in \Omega  \\
 \frac{\partial u}{\partial n}(x)=0,  \ \ \ \ \ \ &  x\in \partial \Omega
\end{array}
\right\}
\end{equation}
where $\varphi_{1}$ is a nontrivial eigenfunction associated to
$\lambda_{1}$ (remember that $\lambda_{1}$ is the first positive
eigenvalue of the eigenvalue problem (\ref{ceq11})). The function
$f(x,u) = \lambda_{1}u + \varphi_{1}(x)$ satisfies {\bf (H)} and
(\ref{ceq10}), but the Fredholm alternative theorem shows that
there is no solution of (\ref{180305}).

\medskip

If, in addition to {\bf (H)} and (\ref{ceq10}), $f$ satisfies a
non-uniform non-resonance condition of the type

\medskip {\bf (h1)} $f_{u}(x,u) \leq \beta (x)$ in $\comega\times
\real$ with $\beta (x) \leq \lambda_{1}$ in $\Omega$ and $\beta
(x) < \lambda_{1}$ in a subset of $\Omega$ of positive measure,

\medskip

\noindent then it has been proved in \cite{mawawi} that
(\ref{ceq9}) has solution. Let us observe that supplementary
condition {\bf (h1)} is given in terms of $\Vert \beta \Vert
_{\infty}.$ In the next result, we provide new supplementary
conditions in terms of $\Vert \beta \Vert _{p},$ where $N/2 <p
\leq \infty$, obtaining a generalization of Theorem 2 in
\cite{mawawi}. For the case $N/2 <p < \infty$ it is clear that the
function $f_u (x,u)$ may cross an arbitrary number of eigenvalues
of the corresponding linear problem, as long as certain
$L_p$-norms are controlled.

In the proof, the basic idea is to combine the results obtained in
the previous section with the Schauder's fixed point theorem.
However we omit this proof here since a more general one will be
given in the case of systems of equations (see Theorem
\ref{nolinealp}).

\begin{theorem}\label{nolineal}
Let $\Omega \subset \real^N$ ($N\geq 2$) be a bounded and regular
domain and $f:\overline{\Omega} \times \real \rightarrow \real $,
$(x,u)\mapsto f(x,u)$, satisfying:
\begin{enumerate}
\item %
$f, f_u$ are
Caratheodory functions and $f(\cdot,0)\in L^{\infty}(\Omega).$%
\item %
There exists a function $\beta \in L^{\infty}(\Omega),$ satisfying
\begin{equation}\label{caratheo}
0 \leq f_{u}(x,u) \leq \beta(x)\ \mbox{in} \ \comega\times \real
\end{equation}
and such that for some $p,$ $N/2 < p \leq \infty,$ we have
$\norm{\beta}_p < \beta_p$
where $\beta_p$ is given by  Theorem \ref{t2partial}.%
\item
 \begin{equation}\label{alphacero} \ \exists
s_0 \in \real \mbox{ s.t. } \displaystyle \int_{\Omega} f(x,s_0) \
dx = 0, \ \mbox{and } f_{u}(x,u(x)) \not\equiv 0, \ \forall u\in
C(\comega)
\end{equation}
\end{enumerate}
Then problem (\ref{ceq9}) has a unique solution.
\end{theorem}

\section{Systems of partial differential equations.}

In this section we deal with elliptic systems of the form

\begin{equation}\label{s1p}
\Delta u(x)+ A(x)u(x) = 0,  \ x\in \Omega,  \ \ \frac{\partial
u(x)}{\partial n}=0, \  x\in \partial \Omega,
\end{equation}
As previously, $\Omega \subset \real^N, \ N \geq 2$ is a bounded
and regular domain and $\dis \frac{\partial}{\partial n}$ is the
outer normal derivative on $\partial \Omega.$ On the other hand,
$A \in \Sigma_{*}$, where $\Sigma_{*}$ is defined as
\afirm{[$\Sigma_{*}$]}{10cm}{The set of real $n\times n$ symmetric
matrix valued function $A(\cdot)$, with continuous element
functions $a_{ij}(x), \ 1 \leq i,j \leq n, \ x \in \comega,$ such
that (\ref{s1p}) has not nontrivial constant solutions and
\begin{equation}\label{04061} \intomega <A(x)k,k>  dx \geq 0, \ \forall \ k \in
\mathbb{R}^n. \end{equation} }  In (\ref{s1p}), $u\in
(H^{1}(\Omega))^n,$ the usual Sobolev space.

$ $

Similar to the ordinary case, to be able to extend to the case of
systems of equations the obtained results in the scalar case, a
key point is the characterization of the best $L_p$ Lyapunov
constant for the scalar given problem, as a minimum of some
especial minimization problems defined in appropriate subsets
$X_p, \ N/2 < p \leq \infty,$ of the Sobolev space $H^1 (\Omega).
$ This is done in the next subsection.  Since $0$ is the first
eigenvalue for Neumann problem, it is necessary to impose an
additional restriction to the definition of the spaces $X_p, \ N/2
<  p \leq \infty$ (see Remark \ref{080311} above). This is not
necessary in the case of Dirichlet one.

\subsection{Lyapunov constants for PDE and constrained minimization problems}

In this subsection we consider $N \geq 2.$ Remember that if
$\frac{N}{2} < p \leq \infty$ and
$$\Sigma = \{ a
\in L^{\infty}(\Omega)\setminus \{0\}: \dis \intomega a(x) \ dx
\geq 0 \ \mbox{and} \ (\ref{s1p}) \ \mbox{has nontrivial
solutions} \ \}
$$ then
$$ \beta_{p} \equiv \min_{a \in \Sigma} \ \Vert a^+ \Vert_p, \ \ N/2 <
p \leq \infty.$$ The next theorem can be proved by using the same
ideas as in the ordinary case (see Lemma 2.5 and Lemma 2.6 in
\cite{camovimia} and Lemma 4.1, Lemma 4.2 and Lemma 4.3 in
\cite{camovijfa}).

\begin{theorem}\label{miap}
If $\frac{N}{2} < p \leq \infty$ is a given number, let us define
the set $X_p$ and the functional $I_p$ as
\begin{equation}\label{eq4p}
\begin{array}{c}
X_{p} = \left \{ v \in H^1(\Omega): \dis \int_{\Omega} \vert v
\vert ^{\frac{2}{p-1}} v = 0 \right \}, \ \ \mbox{if} \ \
\frac{N}{2}
< p < \infty, \\ \\

I_p: X_p \setminus \{ 0 \} \rightarrow \real, \ I_{p}(v) = \dis
\frac{\dis \int_{\Omega} \vert \nabla v\vert^2}{\left ( \dis
\int_{\Omega} \vert v \vert ^{\frac{2p}{p-1}} \right
)^{\frac{p-1}{p}}}, \ \mbox{if} \ \ \frac{N}{2} < p < \infty, \\ \\
X_{\infty} = \{ v \in H^1(\Omega): \ \dis \intomega v = 0 \}, \\
I_\infty: X_\infty \setminus \{ 0 \} \rightarrow \real,\
I_{\infty}(v) = \displaystyle \frac{\dis \intomega \vert \nabla v
\vert ^2}{\displaystyle \intomega v^{2}}
\end{array}
\end{equation}
Then
\begin{equation}\label{eq5p}
\beta_{p} \equiv \min_{X_p \setminus \{ 0 \}} \ I_{p},\
\frac{N}{2} < p \leq \infty,
\end{equation}

\end{theorem}

\begin{remark}\label{rp1} Let us observe that, as in the ordinary case, $\beta_\infty = \lambda_1,$ the first
strictly positive eigenvalue of the Neumann eigenvalue problem in
the domain $\Omega$. Consequently, for PDE it seems difficult to
obtain explicit expressions for $\beta_{p},$ as a function of
$p,\Omega$ and $N,$ at least for general domains. However, for
ODE, it is possible: see (\ref{eqant18}).
\end{remark}

In the next theorem we provide for each $p,$ with $p \in
(N/2,\infty],$ optimal necessary conditions for boundary value
problem (\ref{s1p}) to have only the trivial solution. These
conditions are given in terms of the $L^p$ norm of appropriate
functions $b_{ii}(x), \ 1 \leq i \leq n,$ related to $A(x)$
through the inequality $A(x) \leq B(x), \ \forall\ x \in \comega,$
where $B(x)$ is a diagonal matrix with entries given by
$b_{ii}(x), \ 1 \leq i \leq n.$ In particular, we can use
different $L_{p_i}$ criteria for each $1 \leq i \leq n$ and this
confers a great generality on our results. Since we provide $L^p$
restrictions with $p \in (N/2,\infty),$ the eigenvalues of the
matrix $A(x)$ may cross an arbitrary number of eigenvalues of the
problem
\begin{equation}\label{080311eigenvalue}
\Delta u(x)+ \lambda u(x) = 0,  \ x\in \Omega,  \ \ \frac{\partial
u(x)}{\partial n}=0, \  x\in \partial \Omega,
\end{equation}
 as long as certain $L^p$-norms of the functions $b_{ii}$ are
 controlled. The detailed proof may be seen in \cite{cavijde}.

\begin{theorem}\label{t1p}
Let $A(\cdot) \in \Sigma_{*}$ be such that there exist a diagonal
matrix $B(x)$ with continuous entries $b_{ii}(x),$ and numbers
$p_i \in (N/2,\infty], \ 1 \leq i \leq n,$ which fulfil
\begin{equation}\label{eq7p}
\begin{array}{c}
A(x) \leq B(x), \ \forall \ x  \in \ \comega
\\ \\ \Vert b_{ii}^+
\Vert_{p_i} < \beta_{p_i}, \ 1 \leq i \leq n.
\end{array}
\end{equation}
Then, there exists no nontrivial solution of the vector boundary
value problem (\ref{s1p}).
\end{theorem}

\begin{remark}\label{optimalidadp}
The previous theorem is optimal in the following sense: for any
given positive numbers $\gamma_i, \ 1 \leq i \leq n,$ such that at
least one of them, say $\gamma_j,$ satisfies
\begin{equation}\label{2optimalidad}
\gamma_j > \beta_{p_j}, \ \mbox{for some} \ p_j \in (N/2,\infty],
\end{equation} there exists a diagonal  $n\times n$ matrix
$A(\cdot) \in \Sigma_{*}$ with continuous entries $a_{ii}(x), \ 1
\leq i \leq n,$ satisfying $\Vert a_{ii}^+ \Vert_{p_i} <
\gamma_{i}, \ 1 \leq i \leq n$ and such that the boundary value
problem (\ref{s1p}) has nontrivial solutions.
\end{remark}

 In the next example we consider the case of a system with two
equations. This example is similar to that considered in Remark
\ref{ejemplo240510t} when stability properties for a two
dimensional systems were obtained. Our aim here is to prove that
the unique solution of (\ref{s1p}) is the trivial one, showing the
multiple applications of Lyapunov inequalities.

\begin{example}\label{c2p} Let the matrix $A(x)$ be given by
\begin{equation}\label{eq10}
A(x) = \left (
\begin{array}{cc} a_{11}(x) & a_{12}(x) \\
a_{12}(x) & a_{22}(x) \end{array} \right )
\end{equation}
where \afirm{[{\bf H1}]}{10cm} {
$$
\begin{array}{c} a_{ij} \in C(\comega), \ 1 \leq i,j \leq 2, \\ \\ a_{11}(x) \geq 0, \
a_{22}(x) \geq 0, \ a_{11}(x)a_{22}(x) \geq a_{12}^2 (x), \
\forall \ x \in \comega, \\ \\ \det \ A(x) \neq 0, \ \mbox{for
some}\ x \in \comega. \end{array}
$$ }
In addition, let us assume that there exist $p_1, p_2 \in
(N/2,\infty]$ such that
\begin{equation}\label{eq11}
\Vert a_{11}\Vert_{p_1} < \beta_{p_1}, \ \ \Vert a_{22} +
\displaystyle \frac{a_{12}^2}{\beta_{p_1} - \Vert a_{11}
\Vert_{p_1}} \Vert _{p_2} < \beta_{p_2}.
\end{equation}
Then the unique solution of (\ref{s1p}) is the trivial one.

$ $

To prove the previous affirmation, it is trivial to see that [{\bf
H1}] implies that the eigenvalues of the matrix $A(x)$ are both
nonnegative, which implies that $A(x)$ is positive semi-definite.
Also, since  $\det \ A(x) \neq 0, \ \mbox{for some}\ x \in
\comega,$ (\ref{s1p}) has not nontrivial constant solutions.
Therefore, $A(\cdot) \in \Sigma_{*}.$ Moreover, it is easy to
check that for a given diagonal matrix $B(x),$ with continuous
entries $b_{ii}(x), \ 1 \leq i \leq 2,$ the relation
\begin{equation}\label{eq2306084}
A(x) \leq B(x), \ \forall \ x \in \comega
\end{equation}
is satisfied if and only if $\forall \ x \in \comega,$ we have
\begin{equation}\label{eq12}
\begin{array}{c}
\ b_{11}(x) \geq a_{11}(x), \ b_{22}(x) \geq  a_{22}(x),  \\
(b_{11}(x)- a_{11}(x))(b_{22}(x)- a_{22}(x)) \geq a_{12}^2 (x).
\end{array}
\end{equation}
In our case, if we choose
\begin{equation}\label{eq13}
b_{11}(x) = a_{11}(x) + \gamma, \ b_{22}(x) = a_{22}(x) +
\displaystyle \frac{a_{12}^2 (x)}{\gamma}
\end{equation}
where $\gamma$ is any constant such that
\begin{equation}
\begin{array}{c}
0 < \gamma < \beta_{p_1} - \Vert a_{11} \Vert_{p_1}, \\
\left ( \frac{1}{\gamma} - \frac{1}{\beta_{p_1} - \Vert a_{11}
\Vert _{p_1}} \right ) \Vert a_{12}^2 \Vert_{p_2} < \beta_{p_2} -
\Vert a_{22} + \displaystyle \frac{a_{12}^2}{\beta_{p_1} - \Vert
a_{11} \Vert_{p_1}} \Vert _{p_2}
\end{array}
\end{equation}
then all conditions of Theorem \ref{t1p} are fulfilled and
consequently (\ref{s1p}) has only the trivial solution.
\end{example}

\begin{remark}\label{r2}
The result given in the previous remark may be seen as a
perturbation result in the following sense: let us assume that we
have an uncoupled system of the type
\begin{equation}\label{s2p}
\begin{array}{c}
\Delta u_1  (x) + a_{11}(x)u_1 (x) = 0, \ x \in \Omega; \ \
\frac{\partial u_1 (x)}{\partial n}=0\,  \ x\in \partial \Omega,

\\ \\  \Delta u_2  (x) + a_{22}(x)u_2 (x) = 0, \ x \in \Omega; \ \
\frac{\partial u_2 (x)}{\partial n}=0\,  \ x\in \partial \Omega,
\end{array}
\end{equation}
where \begin{equation}\label{eq14}
\begin{array}{c}
a_{ii} \in C(\comega), \ 1 \leq i \leq 2, \  a_{11}(x) \geq \delta
> 0, \ a_{22}(x) \geq \delta , \ \forall \ x \in \comega.
\\ \\
\exists \ p_1, p_2 \in (N/2,\infty]:\  \Vert a_{11}\Vert_{p_1} <
\beta_{p_1}, \ \ \Vert a_{22}  \Vert _{p_2} < \beta_{p_2}.
\end{array}
\end{equation}
Then it is clear from the scalar results (see Remark \ref{rp1})
that the unique solution of (\ref{s2p}) is the trivial one (see
Corollary 6.1 in \cite{camovijfa}). As it is shown in the previous
remark, we can ensure the permanence of the uniqueness property
(with respect to the existence of solutions) of the coupled system
(\ref{s1p}), for any function $a_{12} \in C(\comega)$ with
$L^\infty-$norm sufficiently small. Here we have considered that
the functions $a_{ii}(x), \ 1 \leq i \leq 2,$ are fixed and that
the uncoupled system is perturbed by the function $a_{12}(x).$ But
it is clear that we may consider, for example, $a_{11}(x),
a_{12}(x)$ fixed and $a_{22}(x)$ as the perturbation.
\end{remark}

\subsection{Resonant nonlinear systems of PDE}

Next we give some new results on the existence and uniqueness of
solutions of nonlinear resonant problems. Similar results can be
proved for ordinary differential systems; in this last case it is
possible to choose the constants $p_i \in [1,\infty], \ 1 \leq i
\leq n $. In particular, next theorem is a generalization, for
systems of equations, of the main result given in \cite{mawawi}
for the Neumann problem. Moreover, it is a generalization (at the
two first eigenvalues) of some results given in \cite{ba} and
\cite{kawa} where the authors take all the constants $p_i =
\infty, \ 1 \leq i \leq n$.

In the proof, the basic idea is to combine the results obtained in
the linear case with Schauder's fixed point theorem.

\bigskip

\begin{theorem}\label{nolinealp}
Let $\Omega \subset \real^N$ ($N\geq 2$) be a bounded and regular
domain and $G : \comega \times \real^n \rightarrow \real, \ (x,u)
\rightarrow G(x,u) $ satisfying:
\begin{enumerate}
\item
\begin{enumerate}
\item $u \rightarrow G(x,u)$ is of class $C^2(\real^n,\real)$ for
every $x \in \comega.$ \item $x \rightarrow G(x,u)$ is continuous
on $\comega$ for every $u \in \real^n.$
\end{enumerate}
\item %
There exist continuous matrix functions $A(\cdot),$ $B(\cdot),$
with $B(x)$ diagonal and with entries $b_{ii}(x),$ and $p_i \in
(N/2,\infty], \ 1 \leq i \leq n,$ such that
\begin{equation}\label{eq15}
\left.
\begin{array}{c}
A(x) \leq G_{uu}(x,u) \leq B(x) \ \mbox{in} \ \comega \times
\real^n, \\ \\ \Vert b_{ii}^+ \Vert_{p_i} < \beta_{p_i}, \ 1 \leq
i \leq n, \\ \\
\intomega <A(x)k,k>  dx > 0, \ \forall \ k \in \mathbb{R}^n
\setminus \{ 0 \}
\end{array}
\right \}
\end{equation}
\end{enumerate}
Then system
\begin{equation}\label{eq17}
\left.
\begin{array}{c}
\Delta u(x)+ G_u (x,u(x)) = 0,  \ x\in \Omega,  \\
\frac{\partial u(x)}{\partial n}=0, \  \ x\in \partial \Omega,
\end{array}
\right \}
\end{equation}
has a unique solution.
\end{theorem}

\begin{proof}  We first prove uniqueness. Let $v$ and $w$ be two
solutions of (\ref{eq17}). Then, the function $u = v-w$ is a
solution of the problem
\begin{equation}\label{eq18}
\Delta u(x) + C(x)u(x)= 0 , \ x \in \Omega, \ \frac{\partial
u}{\partial n}= 0, \ x \in \partial \Omega
\end{equation}
where  $C(x) = \dis \int_{0}^{1} G_{uu}(x,w(x)+\theta u(x)) \ d
\theta$ (see \cite{lang}, p. 103, for the mean value theorem for
the vectorial function $G_u (x,u)$). Hence $A(x)\leq C(x)\leq
B(x)$ and we deduce that $C(x)$ satisfies all the hypotheses of
Theorem \ref{t1p}. Consequently, $u \equiv 0.$
\newline
Next we prove existence. First, we write (\ref{eq17}) in the
equivalent form
\begin{equation}\label{eq19} \left.
\begin{array}{cc}
\Delta u(x) + D(x,u(x))u(x) + G_u(x,0) = 0, & \text{ in } \Omega, \\
 \frac{\partial u}{\partial n}= 0, & \text{ on } \partial \Omega
\end{array}\right\}
\end{equation}
where the  function $D: \comega \times \real^n \rightarrow
\mathcal{M}(\real)$ is defined by $D(x,z) = \dis \int_{0}^{1}
G_{uu}(x,\theta z) \ d \theta.$ Here $\mathcal{M}(\real)$ denotes
the set of real $n \times n$ matrices. Let $X=(C(\comega))^n$ be
with the uniform norm, i.e., if $y(\cdot) = (y^1
(\cdot),\cdots,y^n (\cdot)) \in X, $ then $\Vert y \Vert_X =
\displaystyle \sum_{k=1}^n \Vert y^k (\cdot) \Vert_{\infty}.$
Since
\begin{equation}\label{eq2306085}
A(x) \leq D(x,z) \leq B(x), \ \forall \ (x,z) \in \comega \times
\real^n,
\end{equation}
we can apply Theorem \ref{t1p} in order to have a well defined
operator $T: X \rightarrow X,$ by $Ty = u_{y}$, being $u_{y}$ the
unique solution of the linear problem
\begin{equation}\label{eq20}\left.
\begin{array}{cc}
\Delta u(x)+ D(x,y(x))u(x) + G_u (x,0)= 0, & \text{ in } \Omega, \\
 \frac{\partial u}{\partial n}= 0, & \text{ on } \partial \Omega.
\end{array}\right\}
\end{equation}
Now, we can show that $T$ is completely continuous and that $T(X)$
is bounded (see Theorem 3.4 in \cite{cavijde} for the details).
The Schauder's fixed point theorem provides a fixed point for $T$
which is a solution of (\ref{eq17}).

\end{proof}

\section{Final remarks}

In this section we briefly discuss some questions related to the
previous results. It may be useful to complete the contents of the
paper or for future research.

$ $

{\bf 1.  Lyapunov inequalities and disfocality.}

$ $

When the Neumann (or periodic) boundary value problem is
considered, it is necessary to assume the natural restrictions

\begin{equation}\label{0704111} a\in L^1(0,T)\setminus \{0\}, \ \   \displaystyle
\int_{0}^{T} a(t) \ dt \geq 0. \end{equation}
 Otherwise, we obtain
trivial Lyapunov inequalities (see Remark 4 in \cite{camovimia} ,
Remark 4 in \cite{camovijfa} and Remark \ref{080311} in this
paper).

$ $

Under the restrictions (\ref{0704111}), the relation between, for
example, Neumann boundary conditions, and disfocality arises in a
natural way: if $u \in H^{1}(0,T)$ is any nontrivial solution of
\begin{equation}\label{ceq2} u''(t) + a(t)u(t) =
0, \ t \in (0,T), \ u'(0) = u'(T) = 0,
\end{equation} then $u$ must have a zero $c$ in the interval $(0,T).$
In consequence both problems
$$
 v''(t) + a(t)v(t) = 0, \ t \in (0,c),
\ v'(0) = v(c) = 0, \eqno{{\bf PM(0,c)}}
$$
 and
$$
 v''(t) + a(t)v(t) = 0, \ t \in (c,L),
\ v(c) = v'(L) = 0, \eqno{{\bf PM(c,L)}}
$$
have nontrivial solutions.

$ $

This simple observation can be used to deduce the following
conclusion: if $a \in L^1(0,T)\setminus \{ 0 \}$ with $\int_0^T a
\geq 0$ is any function such that for any $c \in (0,T),$ either
problem {\bf PM(0,c)} or problem {\bf PM(c,L)} has only the
trivial solution, then problem (\ref{ceq2}) has only the trivial
solution. This idea is exploited in \cite{cavidcds} to obtain new
results on the existence and uniqueness of solutions for resonant
problems with Neumann boundary conditions (but it is clear from
the proof that the same idea can be used for other boundary
conditions). For example, by using the $L^\infty$ norm of the
function $a^+$, it can be proved that if
$$
\begin{array}{c}
a \in L^\infty (0,T) \setminus \{ 0 \},\ \displaystyle\int_0^T a
\geq 0 \
\mbox{and} \ \exists \ t_0 \in (0,T): \\ \\
\max \{ t_0^2 \Vert a^+ \Vert_{L^\infty (0,t_0)},  \ (T-t_0)^2
\Vert a^+ \Vert _{L^\infty(t_0,T)} \} \leq \pi^2/4
\end{array}
\eqno{({\bf H})}
$$
and, in addition, either $a^+$ is not the constant $\pi^2 /4t_0^2$
in the interval $[0,t_0]$ or $a^+$ is not the constant $\pi^2
/4(T-t_0)^2$ in the interval $[t_0,T],$ then we obtain that
(\ref{ceq2}) has only the trivial solution (these kinds of
functions $a$ are usually named two step potentials).

Hypothesis ({\bf H}) is optimal in the sense that if $a^+$ is the
constant $\pi^2 /4t_0^2$ in the interval $[0,t_0]$ and $a^+$ is
the constant $\pi^2 /4(T-t_0)^2$ in the interval $[t_0,T],$ then
(\ref{ceq2}) has nontrivial solutions.

If $t_0 = T/2,$ we have the classical result related to the so
called non-uniform non-resonance conditions with respect to the
first two eigenvalues of the Neumann boundary value problem: if
$$ 0 \prec a \prec \frac{\pi^2}{T^2},
$$
then (\ref{ceq2}) has only the trivial solution
(\cite{mawa},\cite{mawa2},\cite{mawawi}). But if for instance $t_0
\in (0,T/2)$ then function $a$ can satisfy $\Vert a^+
\Vert_{L^\infty (0,t_0)} = \pi^2/(4t_0^2)$ (which is a quantity
greater than $ \pi^2/T^2$) as long as $\Vert a^+ \Vert _{L^\infty
(t_0,T)} < \pi^2/4(T-x_0)^2.$

$ $

These ideas are valid not only for $p = \infty,$ but also for
$L_p$ Lyapunov inequalities and disfocality when $1 \leq p <
\infty$ (\cite{cavidcds}) and  also at higher eigenvalues
(\cite{cavijems}).

$ $

{\bf 2. $L_p$-Lyapunov inequalities for scalar equations at higher
eigenvalues, $1 < p < \infty.$ }

$ $

In Section 2, we have obtained the best $L_\infty-L_1$ Lyapunov
constant for the periodic and antiperiodic boundary value problems
at higher eigenvalues. For example, for the periodic problem we
considered the set $\landan$ defined in (\ref{p2806073})  as
\begin{equation}
\landan = \{ a \in \periodica : \lambda_{2n-1} \prec a \
\mbox{and} \ (\ref{p1}) \ \mbox{has nontrivial solutions} \ \}
\end{equation}
and the best $L_\infty-L_1$ Lyapunov constant was defined as
\begin{equation}\label{2pconstante}
\gamma_{1,n} = \inf_{a \in \landan } \ \Vert a  \Vert_{L^1 (0,T)}
\end{equation}

In a similar way, we can obtained very easily the best
$L_\infty-L_\infty$ Lyapunov constant for the periodic problem, at
higher eigenvalues, i.e.,
\begin{equation}
\gamma_{\infty,n} = \inf_{a \in \landan } \ \Vert a
\Vert_{L^\infty (0,T)} = \lambda_{2n+1}
\end{equation}

Remember that $\lambda_{2n-1} = \lambda_{2n} = \frac{4n^2
\pi^2}{T^2}$ and that $\lambda_{2n+1} = \lambda_{2n+2} =
\frac{4(n+1)^2 \pi^2}{T^2}.$

$ $

After this, it is quite natural to consider the same question when
$1 < p < \infty$, that is, try to calculate explicitly
\begin{equation}
\gamma_{p,n} = \inf_{a \in \landan } \ \Vert a \Vert_{L^p (0,T)}
\end{equation}
In this regard, we have developed some ideas but we have no
 a definitive answer yet. The main difficulty is that the associated
minimization problem, defined in (\ref{1603111}) for the case
$p=1,$ is now a convex-concave problem which seems difficult to
study. However, it would be important to have a positive answer to
this question, since in this case, considering $\lim_{p\rightarrow
1^+} \ \gamma_{p,n} $ and $\lim_{p\rightarrow \infty} \
\gamma_{p,n}$ there would be a natural relationship between the
cases $p=1$ and $p= \infty.$

$ $

{\bf 3. Lyapunov inequalities for systems of equations at higher
eigenvalues.}

$ $

In Section 3 we have considered systems of equations and we have
presented some optimal Lyapunov inequalities but always at the
first eigenvalue. For example, for the antiperiodic problem,
$\beta_\infty ^{ant},$ defined in (\ref{eq5}) (see also
(\ref{eqant11})) is equal to $\frac{\pi^2}{T^2},$ the first
eigenvalue of the antiperiodic boundary value problem
(\ref{n2hilleq}) for $a \equiv 0.$

It would be interesting to have $L_1$ Lyapunov inequalities for
systems of equations at higher eigenvalues. The main difficulty is
to obtain a similar characterization, to that given in
(\ref{eqant11}) for the antiperiodic problem, for the constant
$\tilde{\gamma}_{1,n}$ defined in (\ref{ap1009071}).

Finally, but this depends on having been successful or not in the
previous remark, we can try to obtain $L_p$ Lyapunov inequalities
for systems of equations at higher eigenvalues if $1 < p <
\infty.$

$ $

{\bf 4. Explicit value of the Lyapunov constant for PDE and
symmetric domains.}

$ $

In the PDE case, we defined the Lyapunov constant for the Neumann
problem  $\beta_{p,\Omega}$ in (\ref{1603113}) (here we emphasize
the dependence on the domain $\Omega$).

In particular $\beta_{\infty,\Omega} = \lambda_1 (\Omega),$ the
first positive eigenvalue of the eigenvalue problem (\ref{ceq11}).
It is known that except in very special cases, it is not possible
to calculate explicitly $\lambda_1 (\Omega),$ so that the same
will occurs with $\beta_{p,\Omega}.$ However, in Theorem
\ref{miap}, $\beta_{p,\Omega}$ is characterized as the minimum
value of some minimization problems. It is possible that, for
certain types of domains, $\beta_{p,\Omega}$ can be calculated
explicitly. This requires a study of the possible symmetries of
those functions where the minimum is attained. Anyway, this seems
a difficult problem.

$ $

Finally the authors are considering in \cite{caviradial} Lyapunov
inequalities for PDE at radial higher eigenvalues.

\end{document}